%% file: EPOCDG-submitted.tex
\documentclass[final,hidelinks,onefignum,onetabnum]{siamart220329}

\usepackage{mathrsfs}
\usepackage{placeins}
\input{ex_shared}

\makeatletter
\def\refstepcounter@noarg#1{%
  \cref@old@refstepcounter{#1}%
  \cref@constructprefix{#1}{\cref@result}%
  \@ifundefined{cref@#1@alias}%
    {\def\@tempa{#1}}%
    {\def\@tempa{\csname cref@#1@alias\endcsname}}%
  \protected@edef\cref@currentlabel{%
    [\@tempa][\arabic{#1}][\cref@result]%
    \csname p@#1\endcsname\csname the#1\endcsname}%
}
\def\refstepcounter@optarg[#1]#2{%
  \cref@old@refstepcounter{#2}%
  \cref@constructprefix{#2}{\cref@result}%
  \@ifundefined{cref@#1@alias}%
    {\def\@tempa{#1}}%
    {\def\@tempa{\csname cref@#1@alias\endcsname}}%
  \protected@edef\cref@currentlabel{%
    [\@tempa][\arabic{#2}][\cref@result]%
    \csname p@#2\endcsname\csname the#2\endcsname}%
}
\makeatother

\ifpdf
\hypersetup{
	pdftitle={},
	pdfauthor={}
}
\fi

\usepackage{comment}
\usepackage[normalem]{ulem}

\allowdisplaybreaks

\begin{document}
	
	\maketitle
	\setlength{\abovedisplayskip}{2.5pt}   
\setlength{\belowdisplayskip}{2.5pt}

\begin{abstract}
	Central discontinuous Galerkin (CDG) methods offer distinct advantages for hyperbolic conservation laws by evolving solutions on overlapping meshes without requiring complex Riemann solvers or interface numerical fluxes. 
	This design is, however, also the source of a fundamental challenge, as it prevents the direct use of standard summation-by-parts (SBP) flux differencing, making the high-order entropy stability analysis an open problem for CDG methods. Moreover, even within standard DG frameworks where SBP operators are available, existing entropy-stable formulations are mostly semidiscrete, restricted to a single entropy pair, and reliant on positivity- or bound-preserving (BP) assumptions without proof.
	
	In this paper, we present a unified high-order CDG framework (termed EPO-CDG) that simultaneously achieves fully discrete multi-entropy stability, positivity or BP property, and oscillation elimination (OE) near shocks.
	Our analysis reveals that the updated cell averages of standard, arbitrarily high-order CDG schemes with forward Euler time stepping satisfy both a weak version of entropy stability for a general family of multiple convex entropy pairs and a weak version of the BP property under a unified Courant--Friedrichs--Lewy (CFL) condition. 
	Thanks to this key finding, our task reduces to turning the intrinsic weak entropy stability into the desired strong entropy stability. 
	To this end, we construct a unified local convex scaling, termed the weak-to-strong EPO lifting, that enforces multiple discrete entropy inequalities and admissible bounds pointwisely while incorporating a compact OE technique. Coupling this construction with variable-step-size strong-stability-preserving multistep time discretization yields high-order fully discrete schemes with provable multi-entropy stability, without SBP techniques, entropy-conservative numerical fluxes, or Riemann solvers. 
	Extensive numerical tests, ranging from nonconvex scalar conservation laws where single-entropy stability alone can be insufficient to capture the correct entropy solution to challenging compressible flows with strong shocks, demonstrate the effectiveness of the EPO methodology in the CDG framework. 
\end{abstract}

	

\section{Introduction}
For nonlinear hyperbolic conservation laws
\begin{align}\label{equ:conservation_law}
	\bm{u}_t+\nabla\cdot\boldsymbol{F}(\bm{u})=0,
	\qquad
	\boldsymbol{F}=[\bm{f}_1,\ldots,\bm{f}_d],
\end{align}
endowed with strictly convex entropies, discrete entropy inequalities play two important roles. They extend nonlinear stability beyond classical $L^2$ energy bounds and provide the mathematical admissibility criterion needed to identify physically relevant weak solutions. In the scalar setting, enforcing entropy dissipation for all convex entropies uniquely selects the Kru\v{z}kov entropy solution. However, the barrier result of Osher and Tadmor \cite{Osher1984ESchemes,OsherTadmor1988} shows that conservative three-point schemes satisfying this condition are limited to at most first-order accuracy. To overcome this restriction, Tadmor developed a semidiscrete framework based on a single prescribed entropy pair, constructing high-order entropy-stable (ES) schemes by adding suitable dissipation to entropy-conservative (EC) numerical fluxes \cite{tadmor1987numerical,tadmor2003entropy,LeFlochMercierRohde2002,FjordholmMishraTadmor2012TeCNO}. Over the past decade, this framework has been widely developed in the discontinuous Galerkin (DG) methods through the summation-by-parts (SBP) property of Gauss--Lobatto (GL) collocation combined with EC flux differencing or split formulations \cite{CarpenterFisherNielsenFrankel2014,FisherCarpenter2013EntropyStableFD,Gassner2013SkewSymmetricDGSEM}. Subsequent developments include extensions to unstructured simplex meshes \cite{ChenShu2017}, general quadrature rules via hybridized SBP operators and entropy projection \cite{chan2018discretely,ChenShuReview2018}, conservative residual corrections \cite{Abgrall2018GeneralFramework}, entropy-based artificial viscosity \cite{chan2025}, and subcell flux blending \cite{linchen2024}. However, entropy stability in these high-order schemes relies on interface numerical fluxes and typically holds only at the semidiscrete level.

Extending semidiscrete entropy-stable schemes to explicit fully discrete formulations raises additional theoretical difficulties. As Tadmor emphasized \cite[p.~490]{tadmor2003entropy},
\begin{quote}
	\emph{``the entropy stability of explicit schemes hinges on a delicate
		balance between temporal entropy production and spatial entropy
		dissipation.''}
\end{quote}
Standard forward Euler and explicit Runge–Kutta (RK) time discretizations may produce more entropy than is dissipated by the spatial discretization. Existing strategies for achieving fully discrete entropy stability generally modify either the time discretization or the discrete update. For time discretization, implicit time stepping introduces numerical dissipation \cite{LeFlochMercierRohde2002,tadmor2003entropy}, whereas relaxation RK methods adjust the step size or introduce an adaptive scalar parameter to satisfy a discrete entropy inequality \cite{ketcheson2019relaxation,ranocha2020relaxation}. For the discrete update, residual-based corrections \cite{AbgrallOffnerRanocha2022EntropyCorrection} and limiter-based entropy stabilization within the algebraic flux correction framework \cite{kuzmin2022limiter} have been developed to enforce discrete entropy inequalities. Nevertheless, these fully discrete approaches remain mostly restricted to a single entropy pair and still require interface numerical fluxes.

Another important requirement in solving hyperbolic conservation laws is bound preservation (BP), for example maintaining the positivity of density and internal energy in the compressible Euler equations. Such bounds are essential for the entropy pairs considered here, including Harten's entropy family~\cite{harten1983symmetric}, because these entropy functions are defined only for positive density and internal energy \cite{zhang2010positivity,zhang2010maximum}.  In high-order DG discretizations, the canonical BP methodology, pioneered by Zhang and Shu \cite{zhang2010maximum,zhang2010positivity,xu2014bound}, first establishes a weak BP property ensuring that the updated cell averages remain strictly admissible under an explicit Courant--Friedrichs--Lewy (CFL) condition, and then extends this property to the high-order polynomial approximation via a local convex scaling limiter. For complex systems governed by strongly nonlinear constraints, {this methodology was generalized} through the geometric quasilinearization (GQL) approach \cite{wu2023geometric} by transforming nonlinear constraints into equivalent linear ones. In parallel, invariant-domain-preserving (IDP) techniques based on algebraic and convex flux limiting have also been extensively developed \cite{guermond2016invariant,guermond2019invariant,kuzmin2021perspective}. 

Even when a high-order scheme ensures BP and discrete entropy stability, it can still exhibit spurious Gibbs oscillations near discontinuities. Common approaches include nonlinear slope and moment limiters \cite{chavent1989local,cockburn1989tvb,qiu2005runge,zhong2013simple,ShenShe2026MRLimiter}, localized artificial viscosity \cite{hiltebrand2014entropy,huang2020adaptive}, and a posteriori subcell blending or MOOD strategies \cite{clain2011high,Dumbser2014JCP,Dumbser2016JCP}. More recently, jump-based damping was introduced in oscillation-free DG methods \cite{lu2021oscillation,liu2022essentially} and formalized into the scale-invariant oscillation-eliminating (OE) DG (OEDG) framework \cite{peng2023oedg}, which suppresses spurious high-frequency modes through nonintrusive evolution-invariant damping. This mechanism was subsequently integrated with the Zhang--Shu BP limiter via convex oscillation suppression (COS) \cite{CaoEtAl2026COSDG}, demonstrating that local convex scalings can eliminate oscillations and enforce admissibility. However, these convex-scaling limiters can only preserve an established entropy inequality \cite{liu2024entropy,ZuoZhaoLin2024ESPPOFDG} rather than to achieve fully discrete multi-entropy stability. It is challenging to enforce multiple entropy inequalities in a fully discrete scheme while preserving bounds and suppressing spurious oscillations.

Central discontinuous Galerkin (CDG) methods \cite{liu2004central,liu2008stability} offer an alternative to standard DG methods for hyperbolic conservation laws by evolving solutions on staggered, overlapping meshes without requiring interface Riemann solvers or characteristic decompositions. These methods often admit larger CFL numbers while retaining a compact stencil. Over the past decade, BP limiters \cite{LiEtAl2016MPSPPCDG,WuJiangShu2023ProvablyPositiveCDG} and OE techniques \cite{liu2007central,PengWuYuan2025OECDG} have been developed for CDG schemes. {\em Nevertheless, achieving provable high-order entropy stability has remained an open problem, because CDG methods lack the interface numerical fluxes required by standard SBP flux-differencing formulations.}  
Recently, the weak-to-strong stability framework introduced in
\cite{wu2023unified}, which is called EPO, provides a general approach
for unifying entropy stability, positivity, and oscillation suppression. This framework applies
to both finite volume and DG discretizations. Similar to the entropy limiting in EPO, a related limiter was independently
proposed by Liu et al.~\cite{LiuEntropylimiter}.  The entropy limiter is based on the
weak entropy stability of cell averages, which was also established by
Carlier and Renac~\cite{carlier2023invariant} for spectral DG methods. However, extending this approach to CDG discretizations is nontrivial due to their staggered mesh structure. In standard DG formulations, numerical entropy fluxes at cell boundaries enter the balance for the cell average of the entropy, while the SBP structure allows the volume terms to be written in flux-differencing form. In CDG methods, by contrast, the corresponding entropy evolution involves both volume and boundary contributions coupling the primal- and dual-mesh polynomials. This overlapping-mesh coupling precludes a direct application of conventional SBP-based flux decompositions and requires a different treatment of the cell-average entropy inequality.

Our approach is based on the observation that, under a unified CFL condition, the forward Euler update of the CDG cell averages satisfies both weak BP and weak discrete entropy inequalities for a family of strictly convex entropies. We establish this property by decomposing the cell average into a convex combination of subcell GL quadrature values and directional Lax--Friedrichs (LF) intermediate states. These intermediate states are averages of self-similar Riemann solutions over intervals containing the complete wave fans. By Jensen's inequality and the continuous entropy inequalities, these states satisfy the corresponding entropy inequalities for every entropy pair in the chosen family, while the convexity of the invariant admissible domain ensures their admissibility. Thus, weak multi-entropy stability and weak BP follow from a common convex decomposition on staggered meshes. Based on these properties of the updated cell averages, we construct a unified local convex scaling, termed the weak-to-strong EPO lifting, that satisfies multiple discrete entropy inequalities and enforces the prescribed bounds at the quadrature points while incorporating compact primal-dual jump damping for oscillation control. Coupling this spatial lifting with variable-step strong-stability-preserving (SSP) multistep time discretization yields the high-order fully discrete EPO-CDG method.

The proposed framework establishes high-order fully discrete multi-entropy stability without SBP flux differencing or interface Riemann solvers. The resulting compact scheme enforces multiple discrete entropy inequalities simultaneously. This capability helps capture composite waves in nonconvex conservation laws, where a single entropy inequality may fail to select the physically relevant solution.

The primary contributions of this work are summarized as follows:
\begin{itemize}[leftmargin=*]            
\item \emph{Weak multi-entropy stability and its intrinsic connection to weak BP.}
We prove that, under a unified explicit CFL condition, the forward Euler update of the CDG cell averages satisfies both weak BP and weak discrete entropy inequalities for the prescribed family of strictly convex entropies. Both properties follow from the same intrinsic convex decomposition based on averages over Riemann fans. The local scaling uses the admissible cell average and the entropy upper bounds supplied by this decomposition to enforce BP and the prescribed discrete entropy inequalities.
	
	\item \emph{Weak-to-strong EPO lifting under multiple entropy constraints.}
	A unified local convex scaling limiter is developed to extend the weak cell-average ES and BP properties to the high-order polynomial approximation. This construction satisfies multiple discrete entropy inequalities and enforces the prescribed bounds at the quadrature points while incorporating a scale-invariant primal-dual jump damping mechanism to eliminate spurious oscillations without degrading high-order accuracy or spatial compactness.
	
	\item \emph{A high-order, fully discrete EPO-CDG framework and rigorous stability guarantees.}
	By combining the weak-to-strong EPO lifting with an SSP multistep method with variable step sizes, we construct a fully discrete CDG scheme with provable multi-entropy stability and BP. Extensive numerical benchmarks on nonconvex scalar conservation laws, extreme blast waves, and astrophysical jets (with Mach numbers 80 and 2000) validate its designed high-order accuracy, robust shock-capturing capability, and the necessity of enforcing multiple entropy constraints.
\end{itemize}

Throughout this paper, we assume that the admissible state set $\mathcal{G}$ is convex and let $\mathcal{E}$ denote a family of strictly convex entropy functions. 
The remainder of this paper is organized as follows.
Section~\ref{sec:CDG} reviews the semidiscrete CDG formulation on staggered meshes.
In Section~\ref{sec:weak-multi-entropy}, we establish weak multi-entropy stability and its connection to weak BP.
The weak-to-strong EPO lifting and the fully discrete EPO-CDG framework are developed in Section~\ref{sec:EPO}, alongside rigorous stability and admissibility analyses for both forward Euler and high-order SSP multistep time discretizations.
Numerical experiments in Section~\ref{sec:numerical} assess accuracy, entropy stability, and oscillation control, in particular the comparative effects of single and multiple entropy constraints.
Concluding remarks are summarized in Section~\ref{sec:conclusions}.

\section{The CDG method}\label{sec:CDG}

We first review the semidiscrete CDG formulation and specify the notation for the one-dimensional (1D) and two-dimensional (2D) cases. Consider \eqref{equ:conservation_law} with spatial coordinate $\bm x=(x_1,\ldots,x_d)$, written as $x$ in 1D and $(x,y)$ in 2D. Let $\mathcal{T}_h^C$ denote the primal mesh, which partitions the computational domain $\Omega$ into cells. The dual mesh $\mathcal{T}_h^D$ is then obtained by connecting the center of each primal cell to those of its adjacent cells. On these overlapping meshes, we define 
\begin{align*}
	V_h &= \left\{ \bm{v} \in [L^2(\Omega)]^m \;\middle|\; \bm{v}|_K \in [\mathbb{P}^k(K)]^m, \ \forall K \in \mathcal{T}_h^C \right\}, \\
	W_h &= \left\{ \bm{w} \in [L^2(\Omega)]^m \;\middle|\; \bm{w}|_K \in [\mathbb{P}^k(K)]^m, \ \forall K \in \mathcal{T}_h^D \right\},
\end{align*}
Here, $\mathbb P^k(K)$ denotes the space of polynomials of degree at most $k$ on $K$. For any cell $K$, define ${\rm avg}_K(\bm v):=|K|^{-1}\int_K\bm v\,d\bm x$.

On a mesh $\mathcal T$, we define
\begin{equation*}
	\resizebox{\linewidth}{!}{$
			\mathcal G_{\mathcal T}(\bm v,\bm w;\bm\phi_h)
			:=
			\sum_{K\in\mathcal T}
			\Biggl[
			\frac{1}{\tau_{\max}}
			\int_K(\bm w-\bm v)\cdot\bm\phi_h\,d\bm x
			+\int_K\bm F(\bm w):\nabla\bm\phi_h\,d\bm x \\
			-\int_{\partial K}
			\bigl(\bm F(\bm w)\bm n_K\bigr)\cdot\bm\phi_h\,ds
			\Biggr],
		$}
\end{equation*}
where $\bm n_K$ is the outward unit normal, 
$\bm F(\bm w):\nabla\bm\phi_h
=\sum_{\ell=1}^d\bm f_\ell(\bm w)\cdot\partial_{x_\ell}\bm\phi_h$, and
$\bm F(\bm w)\bm n_K
=\sum_{\ell=1}^dn_{K,\ell}\bm f_\ell(\bm w)$.
For general nonlinear fluxes, the integrals in $\mathcal G_{\mathcal T}$ are evaluated via numerical quadrature. Specifically, in the 2D case, the boundary integral over $\partial K$ is approximated by $\langle(\bm F(\bm w)\bm n_K)\cdot\bm\phi_h\rangle_{\partial K}^Q$ using the subcell GL quadrature rule defined in \eqref{eq:epocdg-boundary-face-quadrature}.

For any $\bm{v}_h = (\bm{v}_h^C, \bm{v}_h^D), \bm{\phi}_h = (\bm{\phi}_h^C, \bm{\phi}_h^D) \in V_h \times W_h$, we introduce the $L^2$ inner product on $V_h \times W_h$ as
\begin{equation*}
	\langle \bm{v}_h, \bm{\phi}_h \rangle := \sum_{K \in \mathcal{T}_h^C} \int_K \bm{v}_h^C \cdot \bm{\phi}_h^C \, \mathrm{d}\bm{x} + \sum_{K \in \mathcal{T}_h^D} \int_K \bm{v}_h^D \cdot \bm{\phi}_h^D \, \mathrm{d}\bm{x}.
\end{equation*}
Then, the semidiscrete CDG formulation is given by $\frac{\mathrm{d}\bm{u}_h}{\mathrm{d}t} = \mathcal{L}_h(\bm{u}_h) \in V_h \times W_h$, where 
\begin{align*}
	\langle\mathcal{L}_h(\bm{u}_h),\bm{\phi}_h\rangle
	:=
	\mathcal{G}_{\mathcal{T}_h^C}(\bm{u}_h^C,\bm{u}_h^D;\bm{\phi}_h^C)
	+
	\mathcal{G}_{\mathcal{T}_h^D}(\bm{u}_h^D,\bm{u}_h^C;\bm{\phi}_h^D),
	\quad
	\forall\bm{\phi}_h\in V_h\times W_h.
\end{align*} 
For stability, the maximum allowable time step $\tau_{\max}$ is determined by the CFL condition
$
\tau_{\max} = \frac{C_{\mathrm{CFL}}}{\sum_{\ell=1}^d \alpha_\ell / h_\ell},
$
where $C_{\mathrm{CFL}}$ is the CFL number, $h_\ell$ is the mesh size in the $\ell$th direction, and $\alpha_\ell$ is the corresponding maximum directional wave speed.

\par

In 1D, the primal and dual cells are
$I_j=[x_{j-\frac12},x_{j+\frac12}]$, with center $x_{j}$, and
$I_{j+\frac12}=[x_{j},x_{j+1}]$, respectively, and
$\tau_{\max}=C_{\mathrm{CFL}}h_1/\alpha_1$.
In 2D, the primal and dual cells are
$C_{ij}=[x_{i-\frac12},x_{i+\frac12}]
\times[y_{j-\frac12},y_{j+\frac12}]$, with center $(x_{i},y_{j})$, and
$C_{i+\frac12,j+\frac12}=[x_{i},x_{i+1}]
\times[y_{j},y_{j+1}]$, respectively, and
$\tau_{\max}=C_{\mathrm{CFL}}/
(\alpha_1/h_1+\alpha_2/h_2)$. For the linear advection equation on uniform Cartesian meshes, the semidiscrete CDG method is $L^2$-stable
\cite{liu2008stability} and achieves optimal error estimates \cite{liu2018optimal}.

\paragraph{Quadrature rule}
Each primal or dual cell $K$ is partitioned by the overlapping dual or primal mesh into $2^d$ uniform Cartesian subcells (two subintervals in 1D and four rectangles in 2D), denoted by $\mathcal{S}(K)$. On each subcell, a normalized tensor-product GL quadrature rule is applied. Let $\{(\xi_q, \widehat{\omega}_q)\}_{q=1}^N$ be the $N$-point normalized GL quadrature nodes and weights on the reference interval $[-1, 1]$, satisfying $\sum_{q=1}^N \widehat{\omega}_q = 1$. The number of quadrature points $N$ is chosen such that $2N - 3 \ge k+1$, ensuring that the 1D rule is exact for polynomials of degree up to $k+1$, and the corresponding tensor-product rule integrates all polynomials in $\mathbb{Q}^{k+1}$ exactly on each Cartesian subcell.

\par

For each subcell $S\in\mathcal S(K)$ with directional width $\Delta x_{\ell,S}=h_\ell/2$ (denoted by $\Delta x_S=h/2$ in 1D), let $T_S:[-1,1]^d\to S$ be the affine map from the reference cell to $S$. For a multi-index $\bm q=(q_1,\ldots,q_d)$, set
$\xi_{\bm q}:=(\xi_{q_1},\ldots,\xi_{q_d})$ and $\widehat\omega_{\bm q}:=\prod_{\ell=1}^d\widehat\omega_{q_\ell}.$
In 1D, we denote the quadrature nodes by $x_{S,q}$ and the left and right boundary coordinates by $x_S^L$ and $x_S^R$. In 2D, the nodes are denoted by $(x_{S,\alpha},y_{S,\beta})$, with boundary coordinates $x_S^L,x_S^R$ and $y_S^B,y_S^T$ in the two coordinate directions, respectively.
The quadrature average over $K$ is then defined by the tensor-product GL quadrature averages over all subcells as
\begin{align}
	{\rm avg}_K^Q(v)
	:=
	\frac{1}{2^d}
	\sum_{S\in\mathcal S(K)}
	\sum_{\bm q\in\{1,\ldots,N\}^d}
	\widehat\omega_{\bm q}\,
	v\bigl(T_S(\xi_{\bm q})\bigr). \label{equ:quadrature}
\end{align}
Consequently, $\mathrm{avg}_K^Q(\bm v)=\mathrm{avg}_K(\bm v)$ for every $\bm v\in V_h\cup W_h$.

\begin{remark}
	While tensor-product GL quadrature is a natural choice for decomposing the cell average over the subcells formed by the staggered primal and dual meshes, the resulting decomposition is neither unique nor CFL-optimal. In the standard DG framework, Cui et al.~\cite{CuiDingWu2024OCAD} developed the optimal cell average decomposition theory to improve the CFL condition. Extending the EPO-CDG framework with optimized subcell convex decompositions that relax the CFL restriction will be explored in future work.
\end{remark}

\par

Let $(\eta^{(r)},\mathcal Q^{(r)})$ and $(\eta^{(r)},\bm{\mathcal Q}^{(r)})$ denote the entropy pairs in 1D and 2D, respectively. To formulate the discrete entropy constraints, we denote the cell-average entropy on $K$, approximated via the quadrature rule \eqref{equ:quadrature}, by
\begin{equation*}
	\left\langle \eta^{(r)}(\bm{v}) \right\rangle_K^Q := \mathrm{avg}_K^Q \left( \eta^{(r)}(\bm{v}) \right).
\end{equation*}
In the 1D setting with $K = [x_K^L, x_K^R]$, the boundary term of entropy flux is defined as
\begin{equation*}
	\left\langle \mathcal{Q}^{(r)}(\bm{v}) n_K \right\rangle_{\partial K}^Q
	:=
	-\mathcal{Q}^{(r)}\bigl(\bm{v}(x_K^L)\bigr) + \mathcal{Q}^{(r)}\bigl(\bm{v}(x_K^R)\bigr).
\end{equation*}
Similarly, in the 2D setting with $K = [x_K^L, x_K^R] \times [y_K^B, y_K^T]$, the boundary integral of entropy flux approximated via the subcell quadrature rule is defined by
\begin{equation}\label{eq:epocdg-boundary-face-quadrature}
	\vcenter{\hbox{\resizebox{0.98\hsize}{!}{$
				\begin{aligned}
					\left\langle \bm{\mathcal{Q}}^{(r)}(\bm{v}) \cdot \bm{n}_K \right\rangle_{\partial K}^Q
					&:=
					\frac{h_2}{2} \sum_{\beta=1}^N \widehat{\omega}_\beta \Bigg[
					-\sum_{x_K^L \in \partial S} \bm{\mathcal{Q}}^{(r)}\bigl(\bm{v}(x_K^L, y_{S,\beta})\bigr)
					+ \sum_{x_K^R \in \partial S} \bm{\mathcal{Q}}^{(r)}\bigl(\bm{v}(x_K^R, y_{S,\beta})\bigr)
					\Bigg] \cdot \bm{e}_1
					\\[1mm]
					&\quad + \frac{h_1}{2} \sum_{\alpha=1}^N \widehat{\omega}_\alpha \Bigg[
					-\sum_{y_K^B \in \partial S} \bm{\mathcal{Q}}^{(r)}\bigl(\bm{v}(x_{S,\alpha}, y_K^B)\bigr)
					+ \sum_{y_K^T \in \partial S} \bm{\mathcal{Q}}^{(r)}\bigl(\bm{v}(x_{S,\alpha}, y_K^T)\bigr)
					\Bigg] \cdot \bm{e}_2,
				\end{aligned}
				$}}}
	\raisetag{-2.5em}
\end{equation}
where $\bm{e}_1 = (1, 0)^\top$ and $\bm{e}_2 = (0, 1)^\top$. Here, the notation $\sum_{x_K^L \in \partial S}$ denotes the summation over all subcells $S \in \mathcal{S}(K)$ whose left boundaries lie on \(x=x_K^L\). The remaining boundary sums are defined similarly.

\section{Weak multi-entropy stability and its connection to weak BP}
\label{sec:weak-multi-entropy}

The staggered primal and dual meshes allow an intrinsic convex decomposition of the cell-average update into subcell quadrature point values and LF intermediate states. In this section, we prove that this decomposition yields both weak multi-entropy stability and weak BP for the forward Euler (FE) update. To unify the two analyses, we reformulate the cell-average BP argument \cite{LiEtAl2016MPSPPCDG}, based on the Zhang--Shu framework \cite{zhang2010maximum,zhang2010positivity}, in terms of LF intermediate states \cite{guermond2016invariant,guermond2019invariant}. Let $\theta:=\Delta t/\tau_{\max}\in[0,1]$.
With FE time discretization, the cell average of the CDG solution on each primal cell $K\in\mathcal T_h^C$ satisfies
\begin{equation}\label{eq:app-fe-average}
	\resizebox{0.90\hsize}{!}{$\displaystyle
		{\rm avg}_K(\bm u_h^{C,n+1})
		= (1-\theta){\rm avg}_K(\bm u_h^{C,n})
		+\theta \left[ {\rm avg}_K(\bm u_h^{D,n}) -\frac{\tau_{\max}}{|K|}
		\left\langle \bm F(\bm u_h^{D,n})\cdot\bm n_K \right\rangle_{\partial K}^{Q}\right].
		$}
\end{equation}

We first define weak BP and weak multi-entropy stability for the updated cell averages.
\begin{definition}[Weak BP and weak multi-entropy stability]
	\label{def:epocdg-weak-properties}
	Suppose that $\bm u_h^{C,n}$ and $\bm u_h^{D,n}$ belong to $\mathcal G$ at all subcell {GL} quadrature points. The FE update of the cell averages in CDG schemes is called:
	\begin{enumerate}
		\item \textbf{\textit{Weakly bound-preserving}}, if $\operatorname{avg}_K(\bm u_h^{C,n+1})\in\mathcal G$ for every $K\in\mathcal T_h^C$ and $\operatorname{avg}_K(\bm u_h^{D,n+1})\in\mathcal G$ for every $K\in\mathcal T_h^D$.
		\item \textbf{\textit{Weakly multi-entropy-stable}}, if the discrete entropy inequalities
		\begin{subequations}
			\begin{align}\label{eq:epocdg-weak-entropy-C-unified}
				\eta^{(r)}
				\left(
				\operatorname{avg}_{K}(\bm u_h^{C,n+1})
				\right)
				\le&
				(1-\theta)
				\left\langle
				\eta^{(r)}(\bm u_h^{C,n})
				\right\rangle_{K}^{Q}
				+
				\theta
				\left\langle
				\eta^{(r)}(\bm u_h^{D,n})
				\right\rangle_{K}^{Q}
				\notag\\&\quad
				-
				\frac{\Delta t}{|K|}
				\left\langle
				\bm{\mathcal Q}^{(r)}(\bm u_h^{D,n})\cdot \bm n_K
				\right\rangle_{\partial K}^{Q},
				\qquad K\in\mathcal T_h^C,
			\end{align}
			\begin{align}\label{eq:epocdg-weak-entropy-D-unified}
				\eta^{(r)}
				\left(
				\operatorname{avg}_{K}(\bm u_h^{D,n+1})
				\right)
				\le&
				(1-\theta)
				\left\langle
				\eta^{(r)}(\bm u_h^{D,n})
				\right\rangle_{K}^{Q}
				+
				\theta
				\left\langle
				\eta^{(r)}(\bm u_h^{C,n})
				\right\rangle_{K}^{Q}
				\notag\\&\quad
				-
				\frac{\Delta t}{|K|}
				\left\langle
				\bm{\mathcal Q}^{(r)}(\bm u_h^{C,n})\cdot \bm n_K
				\right\rangle_{\partial K}^{Q},
				\qquad K\in\mathcal T_h^D,
			\end{align}
		\end{subequations}
		hold simultaneously for every $\eta^{(r)}\in\mathcal E$.
	\end{enumerate}
\end{definition}

Let $\lambda_\ell:=2\tau_{\max}\alpha_\ell/h_\ell$ for $\ell=1,\ldots,d$.
For $S\in\mathcal S(K)$, let $\chi_{S,q}^{K}=1$ if $x_{S,q}\in\{x_K^L,x_K^R\}$ and zero otherwise. In 2D, define $\chi_{S,\alpha}^{x,K}=1$ if $x_{S,\alpha}\in\{x_K^L,x_K^R\}$ and $\chi_{S,\beta}^{y,K}=1$ if $y_{S,\beta}\in\{y_K^B,y_K^T\}$, with each indicator zero otherwise.
Using GL quadrature on each subcell as introduced in \Cref{sec:CDG}, together with the directional LF intermediate state
\begin{align*}
	\bm H_\ell(\bm U_L,\bm U_R)
	:=
	\frac12(\bm U_L+\bm U_R)
	-
	\frac{1}{2\alpha_\ell}
	\left(
	\bm f_\ell(\bm U_R)-\bm f_\ell(\bm U_L)
	\right),
\end{align*}
 the FE update of the cell averages  in \eqref{eq:app-fe-average} can be expressed in 1D as
\begin{align}\label{eq:app-fe-average-decomp-1d}
	\mathrm{avg}_K(\bm{u}_h^{C,n+1})
	&=
	(1-\theta)\,\mathrm{avg}_K(\bm{u}_h^{C,n}) + \frac{\theta}{2} \sum_{S\in\mathcal{S}(K)} \sum_{q=1}^N \bigl(\widehat{\omega}_{q} - \lambda_1\chi_{S,q}^{K}\bigr) \bm{u}_h^{D,n}(x_{S,q})
	\notag\\[1mm]
	&\quad + \theta \lambda_1\,\bm{H}_1\big(\bm{u}_h^{D,n}(x_K^L),\, \bm{u}_h^{D,n}(x_K^R)\big),
\end{align}
and in 2D as
\begin{equation}\label{eq:app-fe-average-decomp-2d}
	\resizebox{0.90\hsize}{!}{$
		\begin{aligned}
			{\rm avg}_K(\bm u_h^{C,n+1})
			=&
			(1-\theta){\rm avg}_K(\bm u_h^{C,n})\\
			&\quad+
			\frac{\theta}{4}
			\sum_{S\in\mathcal S(K)}
			\sum_{\alpha,\beta=1}^{N}
			\Bigl(
			\widehat\omega_{\alpha}\widehat\omega_{\beta}
			-
			\lambda_1\widehat\omega_{\beta}\chi_{S,\alpha}^{x,K}
			-
			\lambda_2\widehat\omega_{\alpha}\chi_{S,\beta}^{y,K}
			\Bigr)
			\bm u_h^{D,n}(x_{S,\alpha},y_{S,\beta})
			\\
			&\quad
			+
			\frac{\theta}{2}\lambda_1
			\sum_{\substack{S\in\mathcal S(K)\\x_S^L=x_K^L}}
			\sum_{\beta=1}^{N}
			\widehat\omega_{\beta}\,
			\bm H_1
			\bigl(
			\bm u_h^{D,n}(x_K^L,y_{S,\beta}),
			\bm u_h^{D,n}(x_K^R,y_{S,\beta})
			\bigr)
			\\
			&\quad
			+
			\frac{\theta}{2}\lambda_2
			\sum_{\substack{S\in\mathcal S(K)\\y_S^B=y_K^B}}
			\sum_{\alpha=1}^{N}
			\widehat\omega_{\alpha}\,
			\bm H_2
			\bigl(
			\bm u_h^{D,n}(x_{S,\alpha},y_K^B),
			\bm u_h^{D,n}(x_{S,\alpha},y_K^T)
			\bigr).
		\end{aligned}
		$}
\end{equation}
The same decomposition holds for the cell-average update on dual cells after interchanging the primal and dual meshes.
Under the explicit CFL condition
\begin{align}\label{eq:epocdg-cfl-d}
	\sum_{\ell=1}^d \frac{\tau_{\max}\alpha_\ell}{h_\ell\widehat{\omega}_1/2} \le 1,
\end{align}
all decomposition weights remain nonnegative, thereby yielding admissible convex combinations.

Throughout this paper, we assume that $\mathcal G$ is a convex invariant domain for directional Riemann problems \cite{guermond2016invariant,guermond2019invariant}, and that the corresponding Riemann solutions satisfy the continuous entropy inequalities for all $\eta^{(r)} \in \mathcal E$. Whenever $\alpha_\ell$ is an upper bound on the maximum wave speed in the $\ell$-th direction, the LF intermediate state $\bm H_\ell$ belongs to $\mathcal{G}$ and satisfies multi-entropy inequalities.

\begin{lemma}[Admissibility and entropy inequalities for the LF intermediate state]
	\label{lem:epocdg-H-admissibility}
	For any $\ell\in\{1,\ldots,d\}$ and $\bm U_L,\bm U_R\in\mathcal G$, let $\bm R_\ell(\xi;\bm U_L,\bm U_R)$ denote the self-similar entropy solution of the directional Riemann problem along the $\ell$-th direction. Suppose that $\alpha_\ell>0$ bounds all wave speeds in the sense that $\bm R_\ell(\xi;\bm U_L,\bm U_R)=\bm U_L$ for $\xi<-\alpha_\ell$ and $\bm R_\ell(\xi;\bm U_L,\bm U_R)=\bm U_R$ for $\xi>\alpha_\ell$. Then $\bm H_\ell(\bm U_L,\bm U_R)\in\mathcal G$, and for every $\eta^{(r)}\in\mathcal E$ with associated entropy flux $\bm{\mathcal Q}^{(r)} = (\mathcal Q_1^{(r)},\ldots,\mathcal Q_d^{(r)})^\top$, it holds that
	\begin{equation}\label{eq:epocdg-two-point-entropy-directional}
		\eta^{(r)}\big(\bm H_\ell(\bm U_L,\bm U_R)\big)
		\le \frac{\eta^{(r)}(\bm U_L)+\eta^{(r)}(\bm U_R)}{2}
		- \frac{\mathcal Q_\ell^{(r)}(\bm U_R)-\mathcal Q_\ell^{(r)}(\bm U_L)}{2\alpha_\ell}.
	\end{equation}
\end{lemma}

\begin{proof}
	The admissibility of $\bm H_\ell(\bm U_L,\bm U_R)$ follows Eq.~(2.3) in \cite{carlier2023invariant} and the convexity of $\mathcal G$. Jensen's inequality and the continuous entropy inequalities yields \eqref{eq:epocdg-two-point-entropy-directional}. The details are given in Appendix~\ref{app:weak-ES-BP}.
\end{proof}
\begin{remark}
	In practice, deriving a rigorous analytical upper bound on the exact Riemann wave speeds remains challenging. A conventional choice is to employ the local acoustic speed estimate based on the interface states, $\alpha=\max\{|u_L|+c_L,\,|u_R|+c_R\}$. Otherwise, one can determine an admissible value of $\alpha$ by taking maximum root by solving equation obtained by imposing equality in \eqref{eq:epocdg-two-point-entropy-directional}.
\end{remark}

Combined with \Cref{lem:epocdg-H-admissibility}, the convex decompositions \eqref{eq:app-fe-average-decomp-1d}--\eqref{eq:app-fe-average-decomp-2d} establish weak BP and weak multi-entropy stability for the cell-average update.

\begin{theorem}[Weak multi-entropy stability and weak BP of the FE CDG update]
	\label{thm:epocdg-weak-entropy}
	Let $d \in \{1, 2\}$, and suppose that $\bm u_h^{C,n}$ and $\bm u_h^{D,n}$ belong to $\mathcal G$ at all GL quadrature points on each subcell. For each $\ell=1,\ldots,d$, assume that $\alpha_\ell>0$ bounds the complete fan of every $\ell$-directional Riemann problem arising from the state pairs in both meshes. Then, under the CFL condition \eqref{eq:epocdg-cfl-d}, the FE update of the cell averages is both weakly BP and weakly ES with respect to all entropy pairs associated with $\mathcal E$, in the sense of \cref{def:epocdg-weak-properties}.
	\end{theorem}

\begin{proof}
Under the CFL condition \eqref{eq:epocdg-cfl-d}, equations
\eqref{eq:app-fe-average-decomp-1d}--\eqref{eq:app-fe-average-decomp-2d}
express the updated cell averages as convex combinations of the
quadrature point values of $\bm u_h^{C,n}$ and $\bm u_h^{D,n}$
and the LF intermediate states.
Since all these states lie in $\mathcal G$ by assumption and
\Cref{lem:epocdg-H-admissibility}, weak BP follows directly from
the convexity of $\mathcal G$, as in
\cite{LiEtAl2016MPSPPCDG,zhang2010maximum,zhang2010positivity}.

For the entropy inequalities, it suffices to consider the primal mesh, as the dual mesh follows similarly. Using the convexity of each $\eta^{(r)} \in \mathcal{E}$ and applying Jensen's inequality to \eqref{eq:app-fe-average-decomp-1d} in 1D yields
	\[\resizebox{0.95\hsize}{!}{$\displaystyle
		\begin{aligned}
			\eta^{(r)}({\rm avg}_K(\bm u_h^{C,n+1}))
			\le& (1-\theta)\left\langle\eta^{(r)}(\bm u_h^{C,n})\right\rangle_K^Q
			+\frac{\theta}{2}\sum_{S\in\mathcal S(K)}
			\sum_{q=1}^{N}(\widehat\omega_q-\lambda_1\chi_{S,q}^{K})
			\eta^{(r)}(\bm u_h^{D,n}(x_{S,q}))
			\\
			&+ \theta\lambda_1
			\eta^{(r)}\!\left(\bm H_1(\bm u_h^{D,n}(x_K^L),
			\bm u_h^{D,n}(x_K^R))\right).
		\end{aligned}
		$}\]
	Similarly, applying Jensen's inequality to \eqref{eq:app-fe-average-decomp-2d} in 2D gives
	\[\resizebox{0.95\hsize}{!}{$\displaystyle
		\begin{aligned}
			\eta^{(r)}({\rm avg}_K(\bm u_h^{C,n+1}))
			\le& (1-\theta)\left\langle\eta^{(r)}(\bm u_h^{C,n})\right\rangle_K^Q
			\\&+\frac{\theta}{4}\sum_{S\in\mathcal S(K)}
			\sum_{\alpha,\beta=1}^{N}
			\Bigl(\widehat\omega_{\alpha}\widehat\omega_{\beta}
			-\lambda_1\widehat\omega_{\beta}\chi_{S,\alpha}^{x,K}
			-\lambda_2\widehat\omega_{\alpha}\chi_{S,\beta}^{y,K}\Bigr)
			\eta^{(r)}(\bm u_h^{D,n}(x_{S,\alpha},y_{S,\beta}))
			\\
			&+\frac{\theta}{2}\lambda_1\sum_{\substack{S\in\mathcal S(K)\\x_S^L=x_K^L}}
			\sum_{\beta=1}^{N}\widehat\omega_{\beta}\,
			\eta^{(r)}\!\left(\bm H_1(\bm u_h^{D,n}(x_K^L,y_{S,\beta}),
			\bm u_h^{D,n}(x_K^R,y_{S,\beta}))\right)
			\\
			&+\frac{\theta}{2}\lambda_2\sum_{\substack{S\in\mathcal S(K)\\y_S^B=y_K^B}}
			\sum_{\alpha=1}^{N}\widehat\omega_{\alpha}\,
			\eta^{(r)}\!\left(\bm H_2(\bm u_h^{D,n}(x_{S,\alpha},y_K^B),
			\bm u_h^{D,n}(x_{S,\alpha},y_K^T))\right).
		\end{aligned}
		$}\]       
	In both cases, applying \eqref{eq:epocdg-two-point-entropy-directional} to the LF intermediate states yields the corresponding cell entropy inequalities. The argument on the dual mesh is identical after interchanging the primal and dual meshes.
\end{proof}

\section{Weak-to-strong EPO lifting and fully discrete schemes}\label{sec:EPO}

The EPO-CDG framework integrates two main components: the weak BP and weak multi-entropy stability
properties established in \Cref{sec:weak-multi-entropy} and the local weak-to-strong lifting developed in this section. Specifically, this lifting operator maps the original polynomial to one satisfying the entropy and bound constraints while incorporating OE damping to suppress spurious oscillations. We formulate the EPO lifting for SSP multistep time discretizations, which include the FE method as a special case. Rigorous proofs that the EPO lifting guarantees bound preservation and discrete entropy stability for the polynomial approximation are given in \Cref{subsec:EPO_limiter_analysis}, while the shock-capturing capability of the OE mechanism is demonstrated through numerical experiments in \Cref{sec:numerical}.

\par
\subsection{The fully discrete EPO-CDG method}

Let $\Delta t_n=t^{n+1}-t^n$. At time step $n$, an $s$-step SSP multistep method can be written as
\begin{align}\label{eq:epocdg-vss-sspms}
	\bm u_h^{n+1}
	=
	\sum_{\ell=0}^{s-1}c_{\ell,n}
	\left[
	\widetilde{\bm u}_h^{\,n-\ell}
	+e
	\frac{d_{\ell,n}}{c_{\ell,n}}
	\Delta t_n
	\mathcal L_h\!\left(\widetilde{\bm u}_h^{\,n-\ell}\right)
	\right],
\end{align}
where terms with $c_{\ell,n}=0$ are omitted,
$c_{\ell,n},d_{\ell,n}\ge0$, and
$\sum_{\ell=0}^{s-1}c_{\ell,n}=1$.
In particular, when $s=1$ and $c_{0,n}=d_{0,n}=1$,
the scheme reduces to the FE method. At each prior time level $t^{n-\ell}$ ($0 \le \ell \le s-1$), $\tau_{\max}^{\,n-\ell}$ denotes the maximum admissible FE step size associated with $\widetilde{\bm u}_h^{\,n-\ell}$, as prescribed by the CFL condition \eqref{eq:epocdg-cfl-d}. Following \cite[Theorem~1]{HadjimichaelEtAl2016VSSSSP}, the SSP time-step restriction is given by
\begin{equation}\label{eq:epocdg-ssp-cfl}
	\Delta t_n
	\le
	\min_{\ell:\,d_{\ell,n}>0}
	\frac{c_{\ell,n}}{d_{\ell,n}}\,
	\tau_{\max}^{\,n-\ell}.
\end{equation}

We fix the family of strictly convex entropy functions $ \mathcal{E} $ to be enforced by the EPO lifting, with the associated directional entropy fluxes $\bm{\mathcal{Q}}^{(r)}$ (reducing to the scalar flux $\mathcal{Q}^{(r)}$ in 1D). Using the quadrature notation introduced above, we define the unified discrete boundary integral of the entropy flux for a cell $K$ by 
\begin{align*}
	\mathfrak{F}_K^{(r)}(\bm{v})
	:=
	\begin{cases}
		\left\langle \mathcal{Q}^{(r)}(\bm{v}) n_K \right\rangle_{\partial K}^Q,
		&
		K \subset \mathbb{R},
		\\[0.5em]
		\left\langle \bm{\mathcal{Q}}^{(r)}(\bm{v}) \cdot \bm{n}_K \right\rangle_{\partial K}^Q,
		&
		K \subset \mathbb{R}^2.
	\end{cases}
\end{align*}
At time level $t^n$, for any primal cell $K \in \mathcal{T}_h^C$, we denote by $\widetilde{\bm u}_h^{K,n} := \widetilde{\bm u}_h^{C,n}|_K$ the local approximation on $K$ and by $\widetilde{\bm u}_h^{K^\ast,n} := \widetilde{\bm u}_h^{D,n}|_K$ the dual-mesh approximation restricted to $K$. For a dual cell $K\in\mathcal T_h^D$, the roles of $C$ and $D$ are interchanged.

Accordingly, for the high-order fully discrete scheme, the entropy bound on each cell is defined as the convex combination of those from the FE steps in the multistep decomposition
\begin{equation}\label{eq:epocdg-entropy-bound-ms}
	\begin{aligned}
		E_K^{(r)}(\bm{ u}_h^{n},\ldots,\bm{ u}_h^{n-s+1})
		:=
		\sum_{\ell=0}^{s-1} c_{\ell,n}
		\!\Biggl[
		&\left(1-\theta_{\ell,n}\right)
		\left\langle
		\eta^{(r)}(\bm{\widetilde u}_h^{K,n-\ell})
		\right\rangle_K^Q
		+
		\theta_{\ell,n}
		\left\langle
		\eta^{(r)}(\bm{\widetilde u}_h^{K^\ast,n-\ell})
		\right\rangle_K^Q
		\\
		&-
		\frac{d_{\ell,n}}{c_{\ell,n}}
		\frac{\Delta t_n}{|K|}
		\mathfrak F_K^{(r)}
		\bigl(\bm{\widetilde u}_h^{K^\ast,n-\ell}\bigr)
		\Biggr].
	\end{aligned}
\end{equation}
Here, $\theta_{\ell,n}:=\frac{d_{\ell,n}}{c_{\ell,n}}\frac{\Delta t_n}{\tau_{\max}^{n-\ell}}\le 1$. For the FE time discretization, we write this entropy bound simply as $E_K^{(r)}(\bm u_h^n)$.
\par

\par
\paragraph{\upshape\bfseries The EPO lifting}
When the hypotheses of \Cref{thm:epocdg-weak-entropy} hold at each prior time level $t^{n-\ell}$ and the step-size restriction \eqref{eq:epocdg-ssp-cfl} is satisfied, we have $\operatorname{avg}_K(\bm u_h^{K,n+1})\in\mathcal G$ and $\eta^{(r)}(\operatorname{avg}_K(\bm u_h^{K,n+1}))\le E_K^{(r)}$ for every $\eta^{(r)}\in\mathcal E$, where $E_K^{(r)}$ is the entropy upper bound defined in \eqref{eq:epocdg-entropy-bound-ms}. To enforce bound preservation and the multiple entropy constraints while incorporating OE damping, we apply the local convex scaling
\begin{equation}\label{eq:epocdg-local-scaling-limiter}
	\resizebox{0.90\hsize}{!}{$\displaystyle
		\bm{\widetilde{u}}_h^{K,n+1}(\bm{x})
		:=
		\mathrm{avg}_K(\bm{u}_h^{K,n+1})
		+ 
		\beta_K \left( \bm{u}_h^{K,n+1}(\bm{x}) - \mathrm{avg}_K(\bm{u}_h^{K,n+1}) \right),
		\quad \beta_K \in [0, 1].
		$}
\end{equation}
The composite scaling factor $\beta_K$ is constructed from the BP, ES, and OE scaling factors defined below.

\textbf{BP scaling factor.} The BP scaling factor follows the construction of Zhang and Shu \cite{zhang2010maximum,zhang2010positivity}, as adapted to BP CDG methods in \cite{LiEtAl2016MPSPPCDG}. Specifically, the BP factor $\beta_K^{\mathrm{BP}}\in[0,1]$ is selected from the admissible set
$$
\resizebox{0.98\hsize}{!}{$
	{{\mathcal B_K^{\mathrm{BP}}
	:=}
	\left\{
	\beta \in [0, 1] \;\middle|\;
	\mathrm{avg}_K(\bm{u}_h^{K,n+1}) + \beta \left( \bm{u}_h^{K,n+1}(\bm{x}) - \mathrm{avg}_K(\bm{u}_h^{K,n+1}) \right) \in \mathcal{G}, \; \forall \bm{x} \in \mathcal{N}_K
	\right\}},
	$}
$$
where $\mathcal N_K$ denotes the set of all subcell GL quadrature points in $K$. 
More precisely, we set $\beta_K^{\mathrm{BP}} = \max \mathcal{B}_K^{\mathrm{BP}}$ whenever the maximum is attained, and otherwise choose $\beta_K^{\mathrm{BP}} \in \mathcal{B}_K^{\mathrm{BP}}$ arbitrarily close to $\sup \mathcal{B}_K^{\mathrm{BP}}$ (e.g., when $\mathcal{G}$ is open and the supremum is not attained). The resulting intermediate state is then given by 
$$
\bm{\hat{u}}_h^{K,n+1}(\bm{x}) := \mathrm{avg}_K(\bm{u}_h^{K,n+1}) + \beta_K^{\mathrm{BP}} ( \bm{u}_h^{K,n+1}(\bm{x}) - \mathrm{avg}_K(\bm{u}_h^{K,n+1}) ).
$$

\textbf{ES scaling factor.} For each $\eta^{(r)}\in\mathcal E$, we define the ES scaling factor $\beta_{K,r}^{\mathrm{ent}}\in[0,1]$ relative to the admissible intermediate state $\bm{\hat{u}}_h^{K,n+1}$ by
\begin{align}\label{equ:first_orderESscaling}
	\beta_{K,r}^{\mathrm{ent}}
	:=
	\begin{cases}
		1,
		&
		\left\langle \eta^{(r)}(\bm{\hat{u}}_h^{K,n+1}) \right\rangle_K^Q \le {E_K^{(r)}},
		\\[1em]
		\dfrac{
			{E_K^{(r)}} - \eta^{(r)}\big(\mathrm{avg}_K(\bm{u}_h^{K,n+1})\big)
		}{
			\left\langle \eta^{(r)}(\bm{\hat{u}}_h^{K,n+1}) \right\rangle_K^Q - \eta^{(r)}\big(\mathrm{avg}_K(\bm{u}_h^{K,n+1})\big)
		},
		&
		\text{otherwise}.
	\end{cases}
\end{align}
The overall entropy scaling factor is then defined as $\beta_K^{\mathrm{ent}} := \inf_{\eta^{(r)}\in\mathcal E} \beta_{K,r}^{\mathrm{ent}}.$ By construction, the combined scaling factor $\beta_K^{\mathrm{ent}}\beta_K^{\mathrm{BP}}$ simultaneously guarantees bound preservation in $\mathcal{G}$ and discrete entropy stability for every $\eta^{(r)} \in \mathcal{E}$.

\textbf{OE scaling factor.} Following \cite{PengWuYuan2025OECDG}, we employ an OE scaling factor derived from the overlapping primal--dual mesh representations to preserve the compactness of the CDG framework:
$$
\beta_K^{\mathrm{OE}} := \exp\big( -\Delta t_n\,\delta_K(\bm{u}_h^{n+1}) \big) \in (0, 1],
$$
where $\delta_K(\bm{u}_h) := \max_{1 \le a \le m} \delta_K^{(a)}(\bm{u}_h)$. For each component $a = 1,\ldots,m$, we set $\delta_K^{(a)}(\bm{u}_h) = 0$ if $\sigma_K^{(a)}(\bm{u}_h) = 0$. Otherwise, the componentwise damping rate is defined by
$$
\delta_K^{(a)}(\bm{u}_h) := C_\delta \sum_{e \in \partial K} \frac{\beta_{e,K}|e|}{|K|} \sum_{j=0}^{1} \frac{h_{e,K}^{j}}{\sigma_K^{(a)}(\bm{u}_h)|e|} \int_e \left( \sum_{|\bm{\alpha}|=j} \big| \partial^{\bm{\alpha}}(u_h^{C,(a)} - u_h^{D,(a)}) \big|^2 \right)^{1/2} \mathrm{d}s.
$$
Here $|e| = 1$ for a boundary point of a one-dimensional cell, $h_{e,K} := \sup_{\bm{x} \in K} \operatorname{dist}(\bm{x}, e)$, $\beta_{e,K}$ denotes a local estimate of the maximum wave speed normal to $e$, and $C_\delta > 0$ is a constant depending solely on the polynomial degree (taken as $C_\delta = 1.6, 0.8, 0.6$ for $\mathbb{P}^1, \mathbb{P}^2, \mathbb{P}^3$ elements, respectively). The scale-invariant normalization factor is given by
\begin{equation}
	\sigma_K^{(a)}(\bm{u}_h) :=
	\begin{cases}
		\|u_h^{C,(a)} - \operatorname{avg}_\Omega(u_h^{C,(a)})\|_{L^\infty(\Omega)}, & K \in \mathcal{T}_h^C, \\
		\|u_h^{D,(a)} - \operatorname{avg}_\Omega(u_h^{D,(a)})\|_{L^\infty(\Omega)}, & K \in \mathcal{T}_h^D.
	\end{cases}
\end{equation}

Finally, the unified EPO scaling factor is assembled via
\begin{align}\label{eq:epocdg-beta-def}
	\beta_K := \min\left\{ \beta_K^{\mathrm{OE}},\, \beta_K^{\mathrm{BP}} \beta_K^{\mathrm{ent}} \right\}.
\end{align}
By the convexity of $\mathcal{G}$ and each $\eta^{(r)} \in \mathcal{E}$, the scaled polynomial associated with any $\beta \in [0, \beta_K^{\mathrm{BP}} \beta_K^{\mathrm{ent}}]$ remains admissible in $\mathcal{G}$ and satisfies the discrete multi-entropy inequalities. Consequently, taking the minimum with $\beta_K^{\mathrm{OE}} \in (0, 1]$ preserves rigorous ES and BP while effectively suppressing spurious Gibbs-type oscillations.

\par
\subsection{Properties of the EPO lifting}\label{subsec:EPO_limiter_analysis}
The weak BP and weak multi-entropy stability estimates established in \Cref{sec:weak-multi-entropy} provide the cell-average admissibility and entropy bounds required by the EPO lifting. We now prove that the local convex scaling lifts these properties to the polynomial approximation while preserving cell averages.
\begin{proposition}[Properties of the unified EPO lifting]\label{prop:epocdg-local-scaling}
	Suppose that the admissibility and wave-speed assumptions of \Cref{thm:epocdg-weak-entropy} hold at each time level $t^{n-\ell}$, $\ell=0,\ldots,s-1$, and that the step-size restriction \eqref{eq:epocdg-ssp-cfl} is satisfied. Then the EPO lifting \eqref{eq:epocdg-local-scaling-limiter} with \eqref{eq:epocdg-beta-def} satisfies the following properties on each cell $K$:
	\begin{enumerate}[label=(\arabic*)]
		\item The solution after EPO lifting satisfies
		$$
		\widetilde{\bm u}_h^{{K,}n+1}(\bm x)\in\mathcal G,
		\qquad \forall\,\bm x\in\mathcal N_K.
		$$
		
		\item For every $\eta^{(r)} \in \mathcal{E}$, the lifted solution satisfies the discrete cell entropy inequality
		$$
		\left\langle \eta^{(r)}(\widetilde{\bm u}_h^{K,n+1})\right\rangle_K^Q
		\le E_K^{(r)}(\bm u_h^n,\ldots,\bm u_h^{n-s+1})
		$$
		where $E_K^{(r)}$ is defined in \eqref{eq:epocdg-entropy-bound-ms}. Moreover, under periodic boundary conditions, summing over all primal and dual cells yields the global entropy bound
		\begin{align*}
			&\frac12\sum_{K\in\mathcal T_h^C}|K|
			\left\langle \eta^{(r)}(\widetilde{\bm u}_h^{C,n+1})\right\rangle_K^Q
			+\frac12\sum_{K\in\mathcal T_h^D}|K|
			\left\langle \eta^{(r)}(\widetilde{\bm u}_h^{D,n+1})\right\rangle_K^Q
			\\
			&\quad\le \sum_{\ell=0}^{s-1}c_{\ell,n}\Biggl[
			\frac12\sum_{K\in\mathcal T_h^C}|K|
			\left\langle \eta^{(r)}(\widetilde{\bm u}_h^{C,n-\ell})\right\rangle_K^Q
			+\frac12\sum_{K\in\mathcal T_h^D}|K|
			\left\langle \eta^{(r)}(\widetilde{\bm u}_h^{D,n-\ell})\right\rangle_K^Q
			\Biggr].
		\end{align*}
	\end{enumerate}
\end{proposition}

\begin{proof}
	Statement~(1) follows directly from \eqref{eq:epocdg-local-scaling-limiter}, the definition of $\beta_K^{\mathrm{BP}}$, and the convexity of $\mathcal G$, since $\beta_K\le \beta_K^{\mathrm{BP}}$. 
	
	To prove statement~(2), fix an arbitrary entropy pair $(\eta^{(r)},\bm{\mathcal Q}^{(r)})\in\mathcal E$. We first verify the entropy inequality for the cell averages. For the FE EPO-CDG scheme, \Cref{thm:epocdg-weak-entropy} guarantees that $ \eta^{(r)}\big(\mathrm{avg}_K(\bm{u}_h^{{K,}n+1})\big) \le E_K^{(r)}(\bm{u}_h^n)$. More generally, for an $s$-step SSP multistep method, applying the single-step estimate \eqref{eq:epocdg-weak-entropy-C-unified} to each FE substep in \eqref{eq:epocdg-vss-sspms} and taking the convex combination over a primal cell $K \in \mathcal{T}_h^C$ with the SSP weights $c_{\ell,n}$ yields
	\begin{align*}
		\eta^{(r)}
		\left({\rm avg}_K(\bm u_h^{C,n+1})\right)
		\le
		\sum_{\ell=0}^{s-1} c_{\ell,n}
		\Biggl[
		&\left(1-\theta_{\ell,n}\right)
		\left\langle
		\eta^{(r)}(\bm{\widetilde u}_h^{C,n-\ell})
		\right\rangle_K^Q
		+
		\theta_{\ell,n}
		\left\langle
		\eta^{(r)}(\bm{\widetilde u}_h^{D,n-\ell})
		\right\rangle_K^Q
		\\
		&-
		\frac{d_{\ell,n}}{c_{\ell,n}}
		\frac{\Delta t_n}{|K|}
		\left\langle
		\bm{\mathcal Q}^{(r)}(\bm{\widetilde u}_h^{D,n-\ell})
		\cdot\bm n_K
		\right\rangle_{\partial K}^Q
		\Biggr],
	\end{align*}
	which coincides with the right-hand side of \eqref{eq:epocdg-entropy-bound-ms}. The corresponding estimate on dual cells follows identically upon interchanging the roles of the primal and dual meshes.
	
	Next, with the convention that $\beta_K/\beta_K^{\mathrm{BP}}:=0$ whenever $\beta_K^{\mathrm{BP}}=0$, the convexity of $\eta^{(r)}$, Jensen's inequality and the bound $\beta_K/\beta_K^{\mathrm{BP}} \le \beta_{K,r}^{\mathrm{ent}}$ yield
	\begin{align*}
		\left\langle
		\eta^{(r)}(\widetilde{\bm u}_h^{{K,}n+1})
		\right\rangle_K^Q
		&\le
		\left(
		1-\frac{\beta_K}{\beta_K^{\mathrm{BP}}}
		\right)
		\eta^{(r)}
		\left({\rm avg}_K(\bm u_h^{{K,}n+1})\right)
		+
		\frac{\beta_K}{\beta_K^{\mathrm{BP}}}
		\left\langle
		\eta^{(r)}(\hat{\bm u}_h^{{K,}n+1})
		\right\rangle_K^Q
		\\
		&\le
		\left(1-\beta_{K,r}^{\mathrm{ent}}\right)
		\eta^{(r)}
		\left({\rm avg}_K(\bm u_h^{{K,}n+1})\right)
		+
		\beta_{K,r}^{\mathrm{ent}}
		\left\langle
		\eta^{(r)}(\hat{\bm u}_h^{{K,}n+1})
		\right\rangle_K^Q.
	\end{align*}
	By the definition of $\beta_{K,r}^{\mathrm{ent}}$ in \eqref{equ:first_orderESscaling}, the right-hand side is bounded by $E_K^{(r)}(\bm u_h^n)$ for the FE scheme and by $E_K^{(r)}(\bm u_h^n,\ldots,\bm u_h^{n-s+1})$ for the SSP multistep scheme. This completes the proof of the local entropy bound in statement~(2).
\end{proof}

\section{Numerical tests}\label{sec:numerical}
To assess the high-order accuracy, multi-entropy stability, BP, and OE
performance of the proposed EPO-CDG method, we present a series of 1D
and 2D benchmark problems covering linear advection, nonlinear scalar
conservation laws such as the nonconvex Buckley--Leverett
equation~\cite{liu2024entropy} and the Kurganov--Petrova--Popov (KPP)
problem~\cite{kurganov2007adaptive}, and the compressible Euler equations.
 Unless otherwise specified, we use a relaxation-based SSP Runge--Kutta startup followed by variable-step SSP multistep time discretization of order $k+1$ for accuracy tests with $\mathbb{P}^k$ elements, and of order three for nonsmooth tests. The numerical solutions for the latter problems are visualized in terms of cell averages. Under the CFL conditions established in \cref{thm:epocdg-weak-entropy}, we set the admissible FE step size to $\Delta t_{\mathrm{FE}} = \widehat{\omega}_1 h_1 / (2\alpha_1)$ in 1D and $\Delta t_{\mathrm{FE}} = \widehat{\omega}_1 / [2(\alpha_1/h_1 +\alpha_2/h_2)]$ in 2D, where $\widehat{\omega}_1$ denotes the boundary quadrature weight of the normalized GL rule, and $\alpha_1, \alpha_2$ are the corresponding maximum wave-speed bounds. All simulations are implemented in C/C++ in double precision unless stated otherwise.

To isolate and evaluate the individual roles of BP, ES, and OE, we compare five CDG schemes: EPO-CDG ($\mathrm{ES}+\mathrm{BP}+\mathrm{OE}$), EPCDG ($\mathrm{ES}+\mathrm{BP}$), POCDG ($\mathrm{BP}+\mathrm{OE}$), ECDG ($\mathrm{ES}$ only), and PCDG ($\mathrm{BP}$ only).
Because nonsmooth entropy functions may degrade high-order convergence rates, we enforce discrete entropy stability using the following families of smooth convex entropy pairs.

For scalar conservation laws with Lipschitz continuous fluxes, while the exact Kru\v{z}kov entropies fully characterize weak entropy solutions, their lack of smoothness (Lipschitz continuity only) makes them unsuitable for high-order limiting. 
Therefore, we employ smooth approximate Kru\v{z}kov entropies defined by
\begin{align*}
	\eta_{\alpha,\kappa_s}(u)
	&:=
	(u-\kappa_s)\arctan\bigl(\alpha(u-\kappa_s)\bigr)
	-\frac{1}{2\alpha}\log\left(1+\alpha^2(u-\kappa_s)^2\right),
\end{align*}
where $\alpha>0$ is a smoothing parameter fixed at $20$ in all scalar tests below, and the entropy centers $\{\kappa_s\}_{s=0}^{N_\eta-1}$ in the admissible interval $[u_{\min},u_{\max}]$ are given by $\kappa_0 := (u_{\min}+u_{\max})/2$ for $N_\eta=1$ and $\kappa_s := u_{\min}+\frac{s}{N_\eta-1}(u_{\max}-u_{\min})$ for $N_\eta\ge 2$. 
The associated entropy flux $\mathcal{Q}_{\alpha,\kappa_s}(u)$ is determined up to an additive constant by the compatibility relation $\frac{\mathrm{d}}{\mathrm{d}u}\mathcal{Q}_{\alpha,\kappa_s}(u) = \eta'_{\alpha,\kappa_s}(u) f_1'(u)$.
	
For the compressible Euler equations with state vector $\bm{u} = (\rho, \bm{m}^\top, E)^\top$, we adopt the normalized convex Harten entropy family
	\begin{align*}
		\eta_r(\bm{u}; \bm{u}_0) := -\rho h_r\bigl(S(\bm{u}) - S_0[\bm{u}_0]\bigr),
		\qquad
		\bm{\mathcal{Q}}_r(\bm{u}; \bm{u}_0) := \bm{v} \eta_r(\bm{u}; \bm{u}_0),
	\end{align*}
	where $r < 1/\gamma$, $\bm{v} := \bm{m}/\rho$ is the velocity, $p := (\gamma - 1)(E - \|\bm{m}\|^2/(2\rho))$ is the pressure, and $S(\bm{u}) := \log(p / \rho^\gamma)$ denotes the physical specific entropy. 
	Here, $h_r(z) := (\mathrm{e}^{rz} - 1)/r$ for $r \ne 0$ and \(h_0(z):=z\) otherwise. The initial reference entropy is defined by
	$$
	S_0[\bm u_0]
	:=
	\frac{1}{2|\Omega|}
	\sum_{X\in\{C,D\}}
	\sum_{K\in\mathcal T_h^X}
	|K|\left\langle S(\bm u_{0,h}^K)\right\rangle_K^Q,
	$$
	where $\bm u_{0,h}^K$ denotes the DG initialization on cell $K$.
	Subtracting this reference value from $S(\bm u)$ makes the resulting
	ES scaling factor invariant under shifts.
	
Given an entropy $\eta_r\in\mathcal{E}$, we introduce
\begin{equation}\label{eq:global-discrete-entropy}
	\mathscr{E}_{X,r}^{\,n}
	:=
	\sum_{K\in\mathcal{T}_h^X}
	|K|
	\left\langle
	\eta_r\!\left(\widetilde{\bm{u}}_h^{K,n}\right)
	\right\rangle_K^Q,
	\quad X\in\{C,D\},
	\qquad
	\overline{\mathscr{E}}_r^{\,n}
	:=
	\frac{1}{2}
	\left(
	\mathscr{E}_{C,r}^{\,n}
	+
	\mathscr{E}_{D,r}^{\,n}
	\right).
\end{equation}
For nonperiodic problems, the boundary entropy flux must also be accounted for in the global entropy inequality. 
We therefore define $\overline{\mathscr Q}_{r,\partial\Omega}(t)$ as the average of the discrete outward entropy flux computed on the primal and dual meshes across $\partial\Omega$, which reduces to $\overline{\mathscr Q}_{r,\partial\Omega}(t) = \overline{\mathcal Q}_r(x_R,t) - \overline{\mathcal Q}_r(x_L,t)$ in 1D.
At each step of the variable-step SSP multistep method, we compute
the global discrete entropy residual as
\begin{equation}\label{eq:global-discrete-entropy-residual}
	\mathscr R_r^{\,n+1}
	:=
	\overline{\mathscr E}_r^{\,n+1}
	-
	\sum_{\ell=0}^{s-1}
	c_{\ell,n}\,
	\overline{\mathscr E}_r^{\,n-\ell}
	+
	\Delta t_n
	\sum_{\ell=0}^{s-1}
	d_{\ell,n}\,
	\overline{\mathscr Q}_{r,\partial\Omega}^{\,n-\ell},
\end{equation}
where $c_{\ell,n}$ and $d_{\ell,n}$ are the variable-step SSP coefficients in \eqref{eq:epocdg-vss-sspms}. The global fully discrete entropy inequality holds if and only if $\mathscr R_r^{\,n+1}\le 0$.

\subsection{1D linear advection equation}\label{subsec:num-1d-linear}
We first consider the 1D linear advection equation $u_t+u_x=0$ with periodic boundary conditions. 

\begin{exmp}\label{exmp:linear-smooth}
    For the smooth accuracy test, we take the initial condition $u_0(x)=\sin^2(\pi x)$ on $\Omega=[0,1]$, and compute up to $T=1.1$ with $\Delta t_{\mathrm{FE}}=0.01h_1$. We use the ECDG method with the single quadratic entropy pair $\eta(u)=u^2/2$ and $\mathcal Q(u)=u^2/2$.
    \Cref{tab:epocdg-1d-linear} shows the $L^1$, $L^2$, and $L^\infty$ errors along with the corresponding convergence rates. The results confirm convergence of order $k+1$ across all three discrete norms.
\end{exmp}

\begin{table}[!htp]
    \centering
    \small
    \caption{Errors and convergence rates for the $\mathbb P^k$-based ECDG method for the 1D linear advection problem at $T=1.1$.}
    \label{tab:epocdg-1d-linear}
    \begin{tabular}{c c c c c c c c}
        \toprule
        $k$ & $N$ & $L^1$ error & rate & $L^2$ error & rate & $L^\infty$ error & rate \\
        \midrule
        1 & 40 & 2.21e-03 & -- & 2.48e-03 & -- & 3.97e-03 & -- \\
         & 80 & 2.77e-04 & 2.99 & 3.19e-04 & 2.96 & 5.64e-04 & 2.82 \\
         & 160 & 3.50e-05 & 2.99 & 4.36e-05 & 2.87 & 8.66e-05 & 2.70 \\
         & 320 & 5.56e-06 & 2.66 & 7.00e-06 & 2.64 & 1.48e-05 & 2.54 \\
         & 640 & 1.13e-06 & 2.30 & 1.41e-06 & 2.32 & 3.24e-06 & 2.19 \\
        \midrule
        2 & 40 & 1.20e-05 & -- & 1.85e-05 & -- & 5.94e-05 & -- \\
         & 80 & 6.04e-07 & 4.31 & 9.52e-07 & 4.28 & 4.55e-06 & 3.71 \\
         & 160 & 5.02e-08 & 3.59 & 7.22e-08 & 3.72 & 3.19e-07 & 3.84 \\
         & 320 & 6.28e-09 & 3.00 & 9.03e-09 & 3.00 & 3.99e-08 & 3.00 \\
         & 640 & 7.85e-10 & 3.00 & 1.13e-09 & 3.00 & 4.98e-09 & 3.00 \\
        \midrule
        3 & 40 & 3.21e-08 & -- & 4.38e-08 & -- & 2.19e-07 & -- \\
         & 80 & 2.00e-09 & 4.00 & 2.74e-09 & 4.00 & 1.37e-08 & 4.00 \\
         & 160 & 1.25e-10 & 4.00 & 1.71e-10 & 4.00 & 8.56e-10 & 4.00 \\
         & 320 & 9.51e-12 & 3.72 & 1.19e-11 & 3.85 & 5.40e-11 & 3.99 \\
        \bottomrule
    \end{tabular}
\end{table}

\subsection{1D inviscid Burgers' equation}
\label{subsec:num-1d-burgers}

We next consider the 1D inviscid Burgers' equation $u_t+\left(\frac{u^2}{2}\right)_x=0$ on $x\in[0,2\pi]$ with periodic boundary conditions. 
\begin{exmp}\label{exmp:burgers-smooth}
	For the smooth accuracy test, we apply the $\mathbb{P}^k$ ECDG method with four entropy constraints to the problem whose exact solution is implicitly defined by $u(x,t)=\sin\!\bigl(x-tu(x,t)\bigr)+0.5.$
	The computation is carried out up to $T=0.6$ with $\Delta t_{\mathrm{FE}}=0.01h_1$. The errors and convergence rates are listed in \Cref{tab:epocdg-1d-burgers}, confirming the optimal $(k+1)$th-order accuracy across all norms.
\end{exmp}

\begin{table}[!htp]
	\centering
	\small
	\caption{Errors and convergence rates for the $\mathbb P^k$-based ECDG method solving the 1D inviscid Burgers' equation.}
	\label{tab:epocdg-1d-burgers}
	\begin{tabular}{c c c c c c c c}
		\toprule
		$k$ & $N$ & $L^1$ error & rate & $L^2$ error & rate & $L^\infty$ error & rate \\
		\midrule
		1 & 80 & 7.47e-3 & -- & 5.36e-3 & -- & 9.71e-3 & -- \\
		 & 160 & 1.46e-3 & 2.36 & 1.27e-3 & 2.07 & 2.96e-3 & 1.72 \\
		 & 320 & 3.33e-4 & 2.13 & 3.86e-4 & 1.72 & 9.89e-4 & 1.58 \\
		 & 640 & 5.76e-5 & 2.53 & 7.22e-5 & 2.42 & 2.44e-4 & 2.02 \\
		 & 1280 & 1.01e-5 & 2.51 & 1.27e-5 & 2.51 & 5.77e-5 & 2.08 \\
		\midrule
		2 & 80 & 1.39e-4 & -- & 1.70e-4 & -- & 4.24e-4 & -- \\
		 & 160 & 1.24e-5 & 3.48 & 1.69e-5 & 3.33 & 6.68e-5 & 2.67 \\
		 & 320 & 1.61e-6 & 2.95 & 3.06e-6 & 2.46 & 1.65e-5 & 2.01 \\
		 & 640 & 1.64e-7 & 3.30 & 4.05e-7 & 2.91 & 2.79e-6 & 2.57 \\
		 & 1280 & 1.17e-8 & 3.80 & 3.40e-8 & 3.57 & 3.27e-7 & 3.09 \\
		\midrule
		3 & 80 & 4.23e-6 & -- & 8.22e-6 & -- & 2.60e-5 & -- \\
		 & 160 & 1.19e-7 & 5.15 & 2.53e-7 & 5.02 & 1.53e-6 & 4.08 \\
		 & 320 & 5.81e-9 & 4.36 & 1.48e-8 & 4.10 & 1.00e-7 & 3.94 \\
		 & 640 & 1.69e-10 & 5.10 & 2.54e-10 & 5.86 & 1.67e-9 & 5.91 \\
		\bottomrule
	\end{tabular}
\end{table}

\subsection{1D Buckley--Leverett equation}\label{subsec:num-buckley-leverett}

We consider the 1D Buckley--Leverett equation $u_t+f_1(u)_x=0$ with flux $f_1(u)=\frac{4u^2}{4u^2+(1-u)^2}$ on $\Omega=(-4,4)$, subject to outflow boundary conditions and Riemann initial data $u(x,0)=u_l$ for $x\le0$ and $u(x,0)=u_r$ for $x>0$. 
Because $f_1$ is nonconvex on $[-2,2]$, the entropy solution of the Riemann problem with a single initial discontinuity may consist of multiple elementary waves rather than a single one and is determined by the corresponding convex or concave envelope of the flux. 
Different members of the smooth approximate Kru\v{z}kov entropy family may play distinct roles in constraining the wave structure, and a single entropy pair is therefore generally insufficient to capture the correct entropy solution.

\begin{figure}[!htp]
	\centering
	\captionsetup[subfigure]{skip=1pt}
	\vspace{-10pt}
	\begin{subfigure}[t]{0.33\textwidth}
		\centering
		\vspace{0pt}
		\includegraphics[width=\linewidth]
		{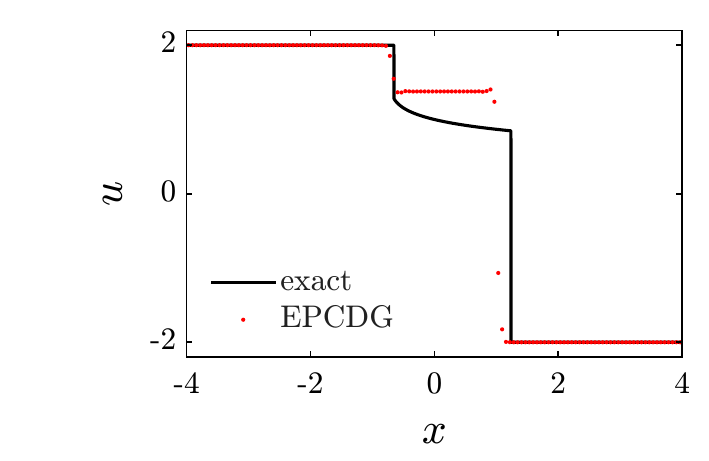}
		\par\vspace{1pt}
		\includegraphics[width=\linewidth]
		{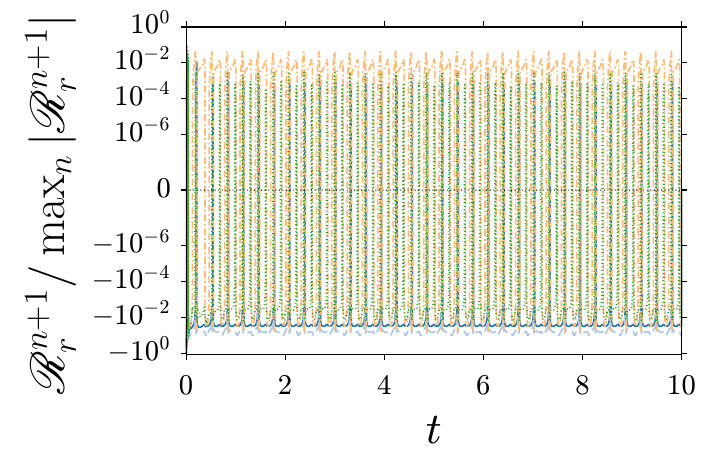}
		\caption{$N_\eta=1$.}
		\label{fig:epcdg-bl-e1}
	\end{subfigure}\hfill
	\begin{subfigure}[t]{0.33\textwidth}
		\centering
		\vspace{0pt}
		\includegraphics[width=\linewidth]
		{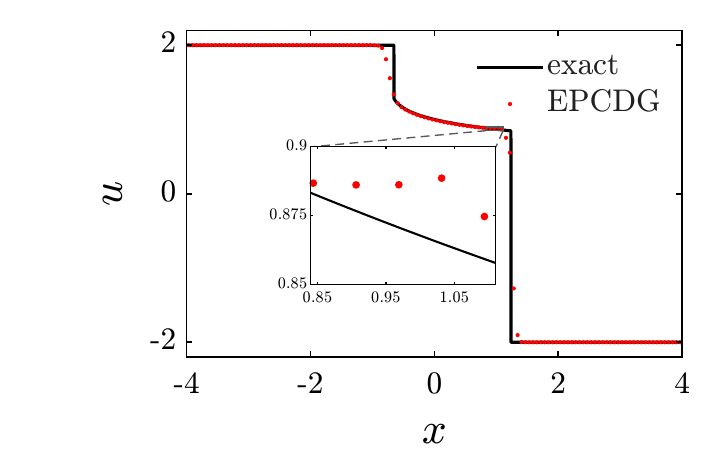}
		\par\vspace{1pt}
		\includegraphics[width=\linewidth]
		{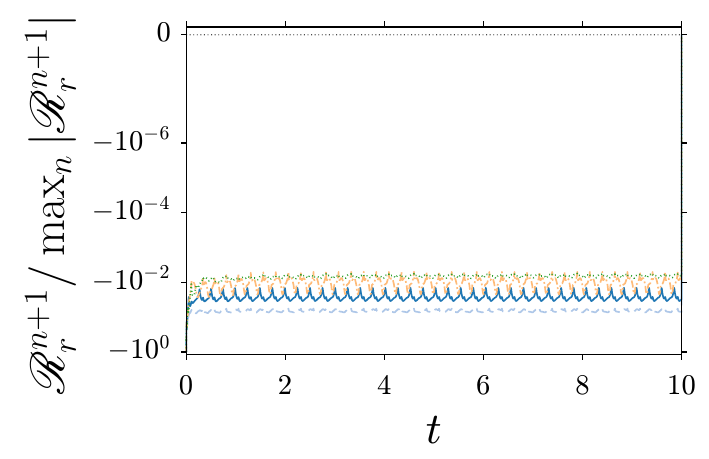}
		\caption{$N_\eta=2$.}
		\label{fig:epcdg-bl-e2}
	\end{subfigure}\hfill
	\begin{subfigure}[t]{0.33\textwidth}
		\centering
		\vspace{0pt}
		\includegraphics[width=\linewidth]
		{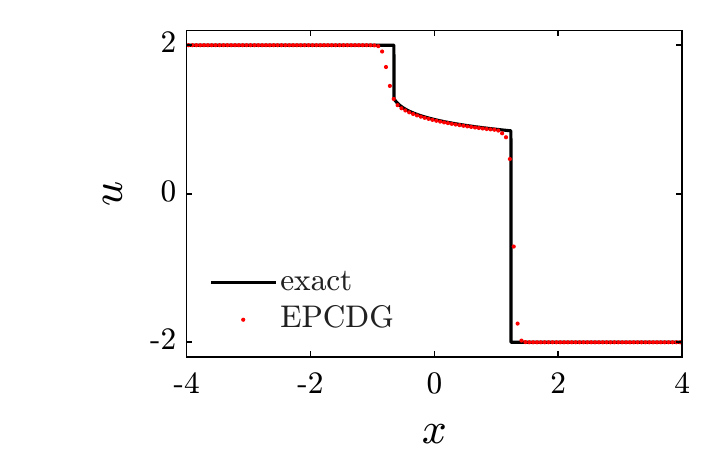}
		\par\vspace{1pt}
		\includegraphics[width=\linewidth]
		{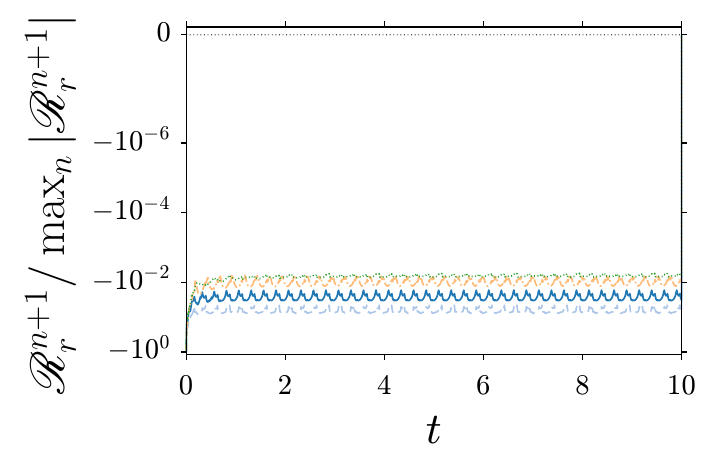}
		\caption{$N_\eta=4$.}
		\label{fig:epcdg-bl-e4}
	\end{subfigure}
	\par\vspace{12pt}
\includegraphics[width=0.90\textwidth]{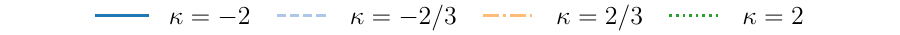}
\caption{EPCDG solutions (top) and normalized global entropy residuals (bottom) for \cref{exmp:bl-2-minus2}.}
	\label{fig:epcdg-bl-entropy-count-comparison}
\end{figure}

\begin{exmp}\label{exmp:bl-2-minus2}
	We consider the Riemann problem with $(u_l,u_r)=(2,-2)$, whose entropy solution exhibits a composite shock-rarefaction-shock structure. We compute up to $T=10$ using the $\mathbb{P}^2$ EPCDG method on $128$ uniform cells. The reference exact solution is sampled at $8001$ uniformly spaced points on $[-4,4]$. We compare the EPCDG solutions obtained with $N_\eta=1$, $2$, and $4$. In all three cases, we evaluate the global discrete entropy residuals $\mathscr R_{\kappa_j}^{\,n+1}$ defined in \eqref{eq:global-discrete-entropy-residual} associated with the four entropy centers $\kappa_j=-2+4j/3$ ($j=0,\ldots,3$), as shown in the bottom row of \cref{fig:epcdg-bl-entropy-count-comparison}.
	
	When only the single entropy centered at $\kappa=0$ is enforced, the numerical solution fails to recover the composite wave structure, and several entropy residuals become positive, most notably for $\kappa_2=2/3$. Enforcing the entropies centered at $\kappa=\pm2$ captures the overall wave structure, with only small local errors near the discontinuities. For $N_\eta=2$, all four global entropy residuals remain nonpositive up to roundoff. With all four entropies enforced, the composite wave structure is well resolved, and all four residuals remain negative at every time step.
\end{exmp}

\subsection{1D compressible Euler equations}\label{subsec:num-1d-euler}

We consider the 1D compressible Euler equations,
$$
\bm U_t+\bm f_1(\bm U)_x=0,
\qquad
\bm U=(\rho,m,E)^{\mathsf T},
\qquad
\bm f_1(\bm U)
=
\left(
m,\frac{m^2}{\rho}+p,\frac{m}{\rho}(E+p)
\right)^{\mathsf T},
$$
where $m=\rho u$. The system is closed by 
$p=(\gamma-1)\left(E-\frac{m^2}{2\rho}\right)$ with $\gamma=1.4.$

\begin{figure}[!htbp]
	\centering
	\makebox[\textwidth][c]{%
		\includegraphics[width=0.50\textwidth,trim=0 0 3bp 0,clip]
		{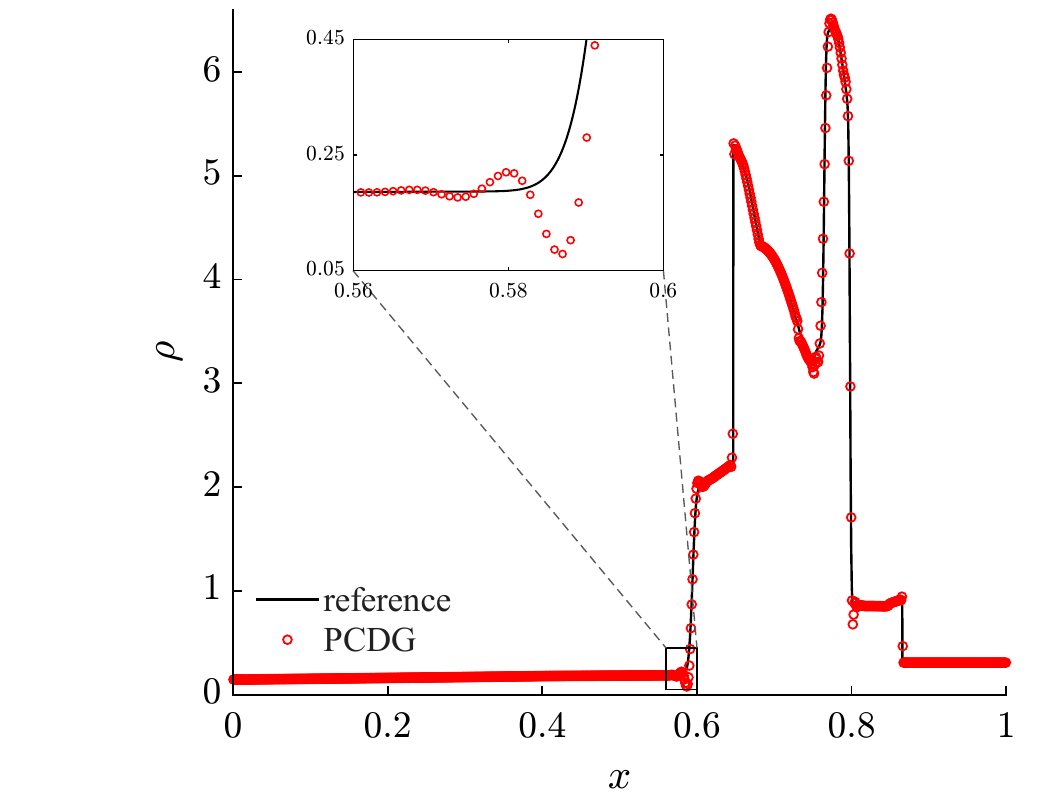}%
		\includegraphics[width=0.50\textwidth,trim=3bp 0 0 0,clip]
		{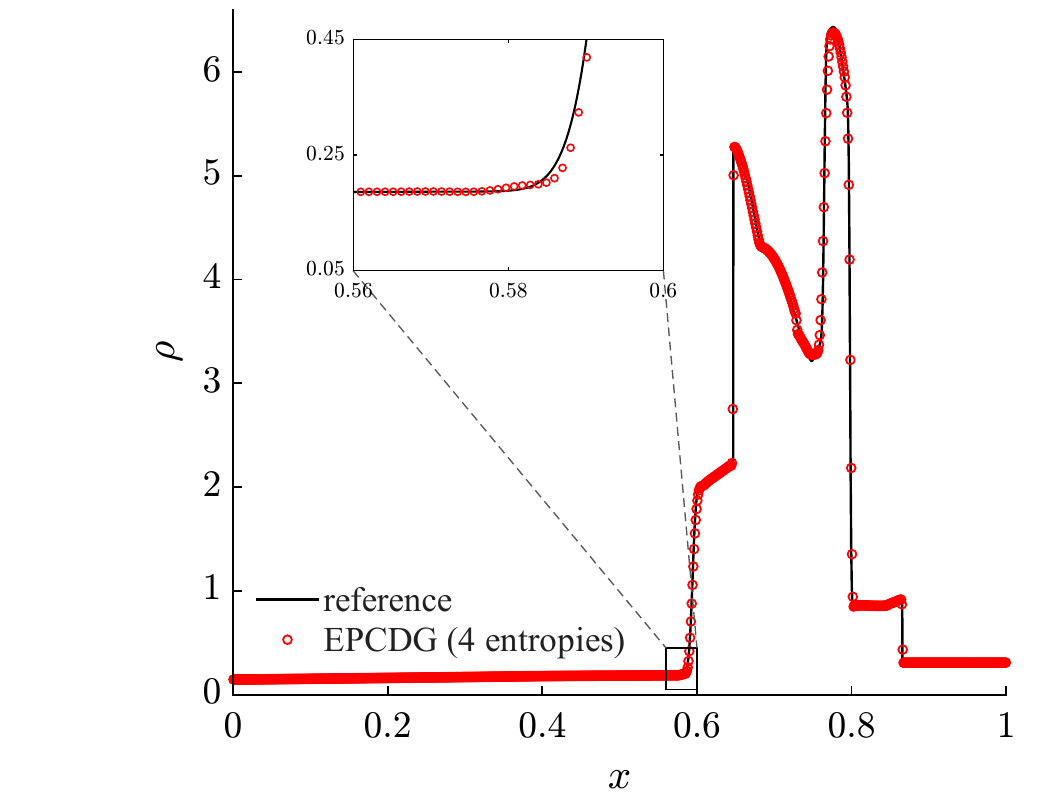}}
	\par\vspace{2pt}
	\makebox[\textwidth][c]{%
		\includegraphics[width=0.50\textwidth,trim=0 0 3bp 0,clip]
		{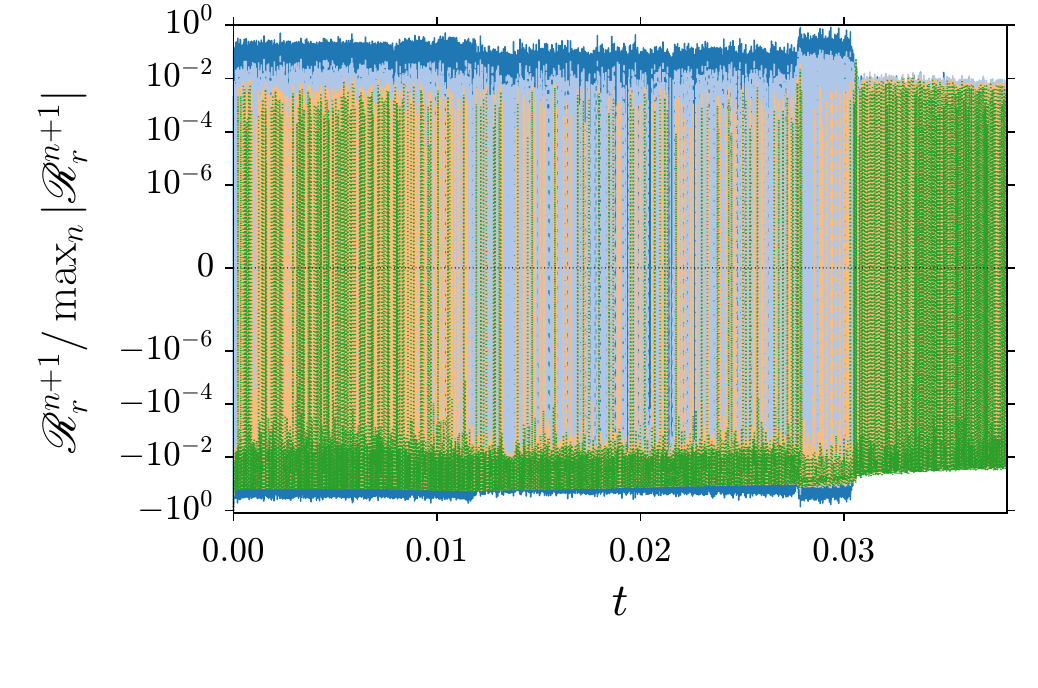}%
		\includegraphics[width=0.50\textwidth,trim=3bp 0 0 0,clip]
		{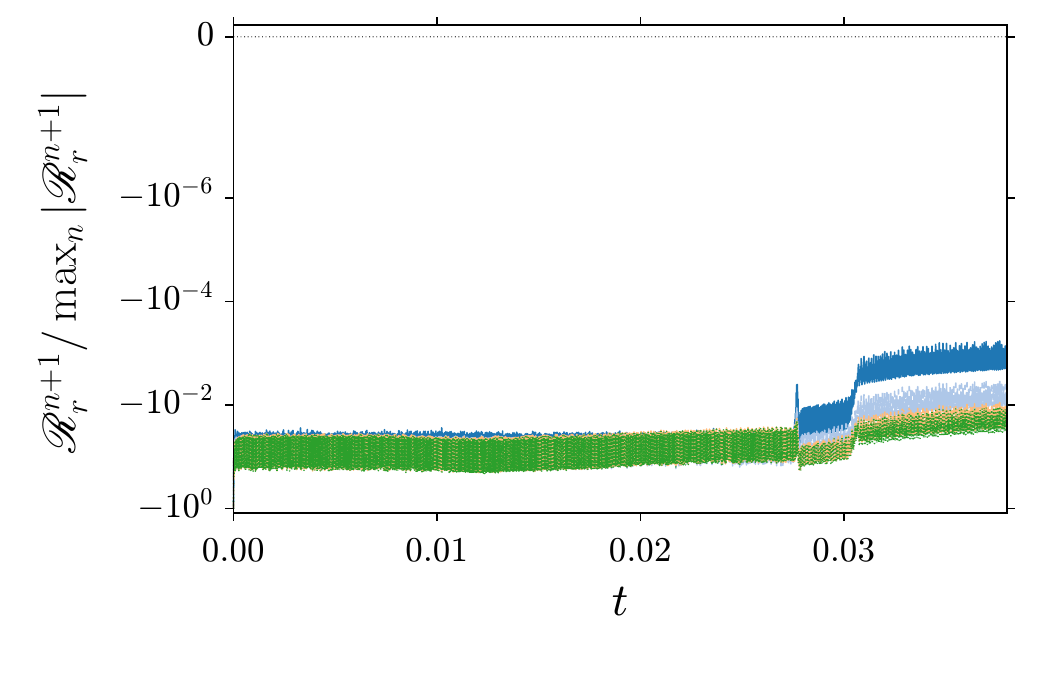}}
	\vskip-10pt
	\par\vspace{1pt}
\includegraphics[width=0.90\textwidth]{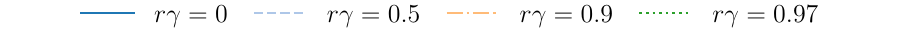}
\caption{Density (top) and normalized global entropy residuals (bottom) for \cref{exmp:euler-two-blast}. }
\label{fig:two-blast-solution-entropy-residual}
\end{figure}

\begin{exmp}[Two blast waves \cite{woodward1984numerical}]\label{exmp:euler-two-blast}
	On the domain $\Omega = [0,1]$, we initialize $\rho(x,0) = 1$, $u(x,0) = 0$, and let $p(x,0)$ take the values $1000$, $0.01$, and $100$ on $0\le x\le0.1$, $0.1<x\le0.9$, and $0.9<x\le1$, respectively, with reflecting boundary conditions imposed at both ends. We compute up to $T=0.038$ using $\mathbb{P}^2$ elements on $960$ uniform cells, comparing the PCDG method with the four-entropy EPCDG method ($r \in \{0,\, 0.50/\gamma,\, 0.90/\gamma,\, 0.97/\gamma\}$) against a reference solution computed by the finite-volume local Lax--Friedrichs (FVLLF) scheme on $300{,}000$ uniform cells. We also evaluate the global discrete entropy residuals \eqref{eq:global-discrete-entropy-residual}. As shown in \cref{fig:two-blast-solution-entropy-residual}, the PCDG solution develops spurious oscillations near the discontinuities and yields positive entropy residuals at several time steps, violating the fully discrete entropy inequalities. In contrast, the four-entropy EPCDG solution agrees closely with the reference solution, while all four global entropy residuals remain nonpositive at every time step.
\end{exmp}

\begin{figure}[!htbp]
	\centering
	\makebox[\textwidth][c]{%
		\includegraphics[width=0.49\textwidth]
		{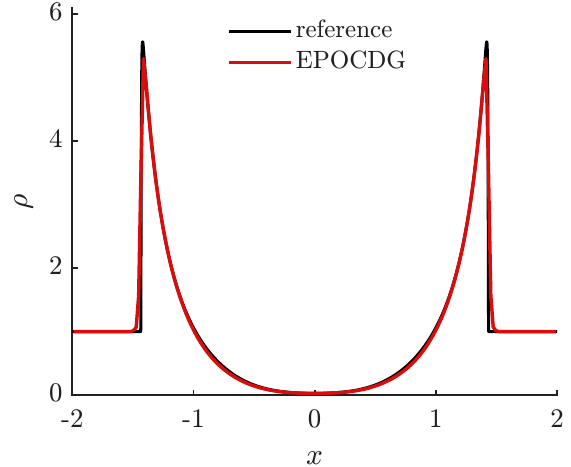}%
		\hspace{0.01\textwidth}%
		\raisebox{0.002710\textwidth}{%
		\includegraphics[width=0.572256\textwidth]
		{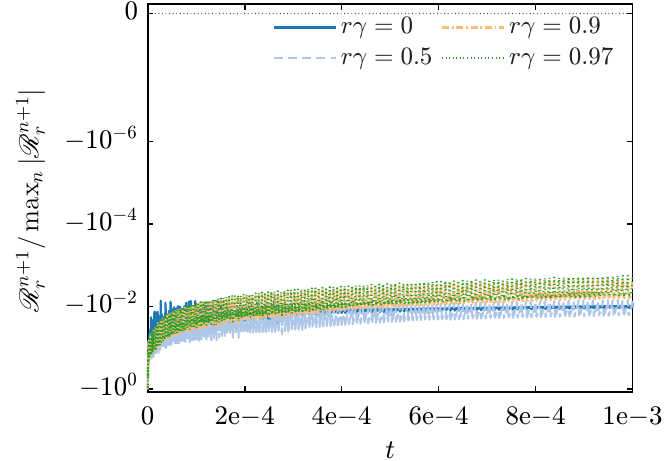}}}
	\caption{Density (left) and normalized global entropy residuals (right) for the four-entropy $\mathbb{P}^2$ EPO-CDG solution in \cref{exmp:euler-sedov}.}
	\label{fig:epocdg-euler-sedov}
\end{figure}

\begin{exmp}[Sedov blast wave \cite{zhang2010positivity}]\label{exmp:euler-sedov}
On $\Omega=[-2,2]$, we set $\rho(x,0)=1$ and $m(x,0)=0$, with $E(x,0)=E_0/h_1$ for $|x|\le h_1/2$ and $E(x,0)=10^{-12}$ otherwise. 
Here, we use $h_1$ for the mesh size, set $E_0=3.2\times10^6$, and impose outflow boundary conditions at $x=\pm2$. 
We then compute the solution up to $T=10^{-3}$ using $\mathbb P^2$ elements on $201$ uniform cells. 
As shown in \Cref{fig:epocdg-euler-sedov}, the EPO-CDG solution using the four entropy pairs with $r\in\{0,0.50/\gamma,0.90/\gamma,0.97/\gamma\}$ captures the two shocks and agrees closely with the FVLLF reference solution sampled at $8001$ points on $[-2,2]$. Moreover, all four global discrete entropy residuals remain nonpositive throughout the computation.
\end{exmp}

\begin{figure}[!htbp]
	\centering
	\makebox[\textwidth][c]{%
		\includegraphics[width=0.49\textwidth]
		{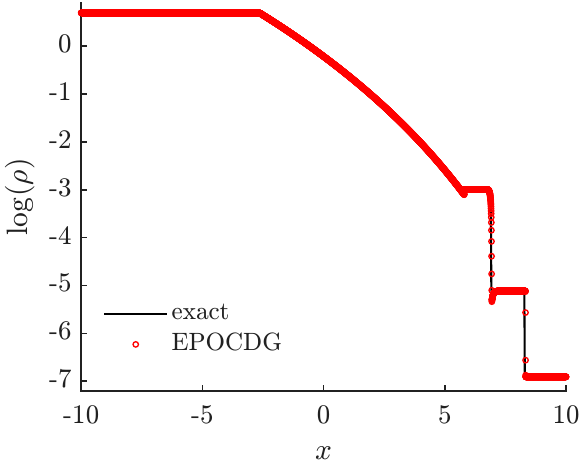}%
		\hspace{0.01\textwidth}%
		\raisebox{0.002891\textwidth}{%
		\includegraphics[width=0.558000\textwidth]
		{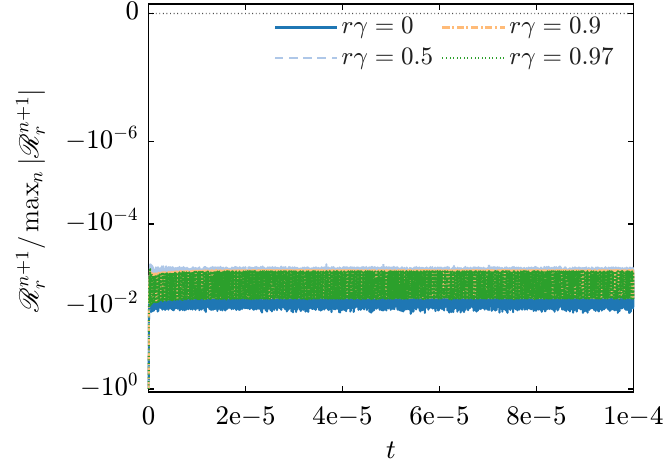}}}
\caption{Logarithmic density (left) and normalized global entropy residuals (right) for the four-entropy $\mathbb{P}^2$ EPO-CDG solution in \cref{exmp:euler-leblanc}.}
	\label{fig:epocdg-euler-leblanc}
\end{figure}

\begin{exmp}[LeBlanc shock tube \cite{xu2022third}]\label{exmp:euler-leblanc}
	On $\Omega=[-10,10]$, we take the initial data $(\rho,u,p)(x,0)=(2,0,10^9)$ for $x<0$ and $(\rho,u,p)(x,0)=(10^{-3},0,1)$ for $x\ge0$. This is an extreme Riemann problem with large jumps in $\rho$ and $p$. Outflow boundary conditions are imposed. The solution is computed up to $T=10^{-4}$ using the four-entropy $\mathbb{P}^2$ EPO-CDG method ($r \in \{0,\, 0.50/\gamma,\, 0.90/\gamma,\, 0.97/\gamma\}$) on 6400 uniform cells. The exact Riemann solution is sampled at $1000$ uniformly spaced points on $[-10,10]$. As shown in \Cref{fig:epocdg-euler-leblanc}, the proposed method sharply resolves the extreme density variation, agrees closely with the exact solution, and keeps all four global discrete entropy residuals nonpositive at every time steps.
\end{exmp}

\begin{table}[!htbp]
	\centering
	\small
	\caption{Errors and convergence rates for the 2D linear advection problem at $T=1.1$, computed using the $\mathbb{P}^k$-based EPO-CDG method.}
	\label{tab:epocdg-2d-linear-accuracy}
	\setlength{\tabcolsep}{7pt}
	\renewcommand{\arraystretch}{1.04}
	\begin{tabular}{@{}ccrcrcrc@{}}
		\toprule
		$k$ & $N$ &
		$L^1$ error & rate &
		$L^2$ error & rate &
		$L^\infty$ error & rate \\
		\midrule
		1 & 16 & $2.00\mathrm{e}{-1}$ & -- & $2.22\mathrm{e}{-1}$ & -- & $3.15\mathrm{e}{-1}$ & -- \\
		& 32 & $3.70\mathrm{e}{-2}$ & 2.43 & $4.15\mathrm{e}{-2}$ & 2.42 & $6.10\mathrm{e}{-2}$ & 2.37 \\
		& 64 & $4.68\mathrm{e}{-3}$ & 2.98 & $5.53\mathrm{e}{-3}$ & 2.91 & $9.11\mathrm{e}{-3}$ & 2.74 \\
		& 128 & $5.68\mathrm{e}{-4}$ & 3.04 & $7.07\mathrm{e}{-4}$ & 2.97 & $1.27\mathrm{e}{-3}$ & 2.84 \\
		& 256 & $7.21\mathrm{e}{-5}$ & 2.98 & $9.13\mathrm{e}{-5}$ & 2.95 & $2.03\mathrm{e}{-4}$ & 2.65 \\
		\midrule
		2 & 16 & $5.08\mathrm{e}{-3}$ & -- & $5.47\mathrm{e}{-3}$ & -- & $7.12\mathrm{e}{-3}$ & -- \\
		& 32 & $3.35\mathrm{e}{-4}$ & 3.92 & $3.52\mathrm{e}{-4}$ & 3.96 & $5.37\mathrm{e}{-4}$ & 3.73 \\
		& 64 & $2.07\mathrm{e}{-5}$ & 4.01 & $2.36\mathrm{e}{-5}$ & 3.90 & $4.35\mathrm{e}{-5}$ & 3.63 \\
		& 128 & $1.33\mathrm{e}{-6}$ & 3.96 & $1.54\mathrm{e}{-6}$ & 3.94 & $3.27\mathrm{e}{-6}$ & 3.73 \\
		& 256 & $9.50\mathrm{e}{-8}$ & 3.81 & $1.18\mathrm{e}{-7}$ & 3.70 & $2.94\mathrm{e}{-7}$ & 3.47 \\
		\midrule
		3 & 16 & $1.32\mathrm{e}{-4}$ & -- & $1.56\mathrm{e}{-4}$ & -- & $3.23\mathrm{e}{-4}$ & -- \\
		& 32 & $4.89\mathrm{e}{-6}$ & 4.76 & $5.66\mathrm{e}{-6}$ & 4.78 & $1.22\mathrm{e}{-5}$ & 4.73 \\
		& 64 & $1.70\mathrm{e}{-7}$ & 4.85 & $1.97\mathrm{e}{-7}$ & 4.85 & $5.04\mathrm{e}{-7}$ & 4.60 \\
		& 128 & $6.14\mathrm{e}{-9}$ & 4.79 & $7.28\mathrm{e}{-9}$ & 4.76 & $2.38\mathrm{e}{-8}$ & 4.41 \\
		& 256 & $2.67\mathrm{e}{-10}$ & 4.52 & $3.23\mathrm{e}{-10}$ & 4.49 & $1.25\mathrm{e}{-9}$ & 4.25 \\
		\bottomrule
	\end{tabular}
\end{table}

\subsection{2D linear advection equation}\label{subsec:num-2d-linear}

We consider the 2D linear advection equation $u_t+u_x+u_y=0$ on Cartesian overlapping meshes with periodic boundary conditions.

\begin{exmp}\label{exmp:2d-linear-smooth}
	On $\Omega=[0,1]^2$, we use the exact solution $u(x,y,t)=\sin^2\!\bigl(\pi(x+y-2t)\bigr)$. For this smooth problem, ES is imposed with the quadratic entropy pair $\eta(u)=\frac{u^2}{2}$ and $\bm{\mathcal Q}(u)=\left(\frac{u^2}{2},\frac{u^2}{2}\right)$. We compute to $T=1.1$ with $\Delta t_{\mathrm{FE}}=0.005h$ using $\mathbb{P}^1$, $\mathbb{P}^2$, and $\mathbb{P}^3$ approximations on $N\times N$ uniform meshes, where $h=1/N$. The $L^1$, $L^2$, and $L^\infty$ errors and observed convergence rates are given in \Cref{tab:epocdg-2d-linear-accuracy}. The rates confirm the high-order accuracy of the EPO-CDG method for this smooth problem.
\end{exmp}

\subsection{2D KPP problem}\label{subsec:num-2d-kpp}

We consider the 2D KPP equation $u_t+(\sin u)_x+(\cos u)_y=0$ with nonconvex flux \cite{kurganov2007adaptive}.

\begin{figure}[!htbp]
	\centering
	\includegraphics[width=0.98\textwidth]{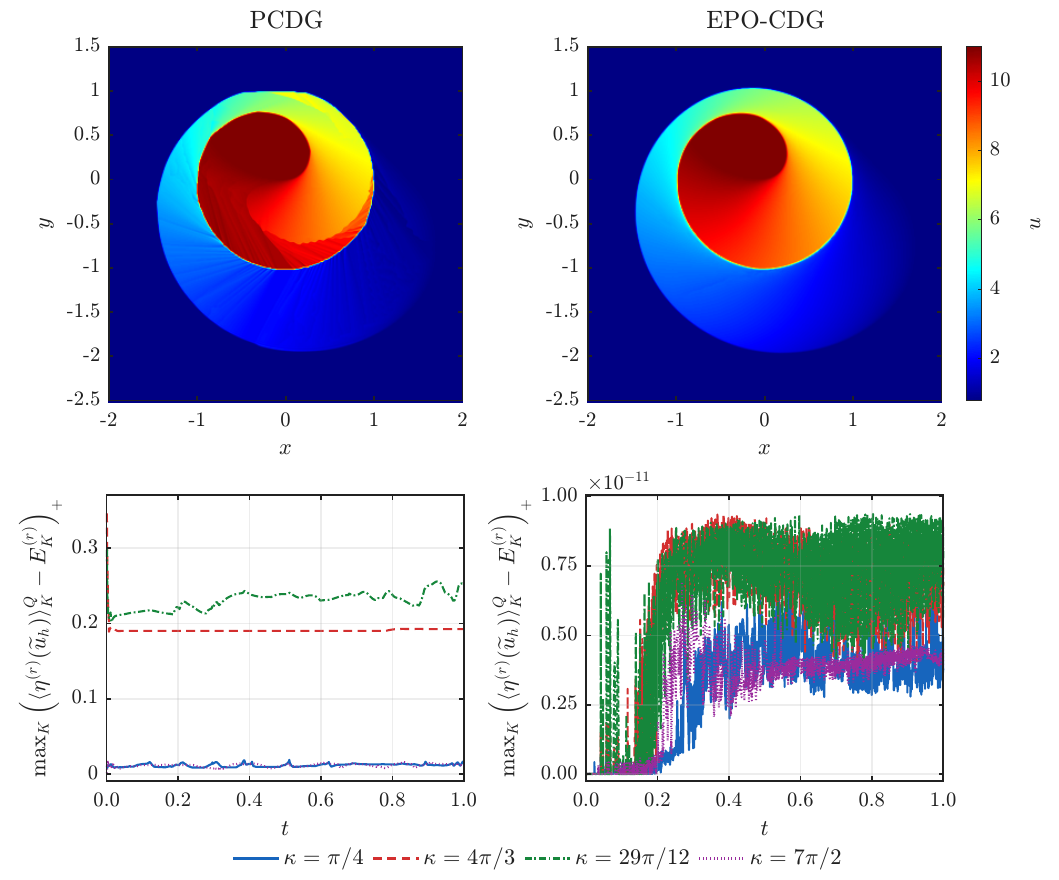}
	\caption{Cell averages at $T=1$ (top) and maximum positive local discrete entropy residuals on the primal mesh (bottom) for \cref{exmp:2d-kpp}, using PCDG (left) and EPO-CDG (right). }
	\label{fig:epocdg-2d-kpp-four-settings}
\end{figure}

\begin{exmp}\label{exmp:2d-kpp}
	On the domain $\Omega=[-2,2]\times[-2.5,1.5]$, we prescribe the discontinuous initial condition $u(x,y,0)=3.5\pi$ for $x^2+y^2<1$ and $u(x,y,0)=0.25\pi$ for $x^2+y^2\ge1$ with periodic boundary conditions imposed on all four sides. We compute up to $T=1$ using $\mathbb{P}^2$ elements on a $320\times320$ uniform mesh. For EPO-CDG, discrete entropy stability is enforced using four smooth approximate Kru\v{z}kov entropies centered at $\kappa \in \{\pi/4,\, 4\pi/3,\, 29\pi/12,\, 7\pi/2\}$.
	We compare the $\mathbb{P}^2$ PCDG and EPO-CDG methods in \cref{fig:epocdg-2d-kpp-four-settings}, where the top row shows the cell-average solutions at $T = 1$ and the bottom row displays the time evolution of the maximum positive local discrete entropy residual on the primal mesh for each prescribed entropy $\eta^{(r)}$,
$$
\max_{K\in\mathcal T_h^C}\max\left\{\left\langle\eta^{(r)}(\widetilde u_h^{C,n+1})\right\rangle_K^Q-E_K^{(r)},\,0\right\},
$$
	with $E_K^{(r)}$ denotes the cell entropy upper bound in \eqref{eq:epocdg-entropy-bound-ms}. 
	While PCDG exhibits noticeable violations of the discrete entropy inequalities, the maximum positive local entropy residuals for EPO-CDG remain within machine roundoff throughout the simulation.
\end{exmp}


\subsection{2D compressible Euler equations}\label{subsec:num-2d-euler}

We next consider several 2D benchmark problems for the compressible Euler equations with an ideal-gas equation of state. Unless otherwise specified, the ratio of specific heats is set to $\gamma=1.4$. 

\begin{exmp}\label{exmp:2d-euler-density-wave}
	On $\Omega=[0,2]^2$, we consider the smooth stationary density wave with $\rho(x,y,t)=1+0.2\sin\!\bigl(\pi(x-y)\bigr)$ with $u=v=p=1$ and periodic boundary conditions.
	The $\mathbb{P}^1$, $\mathbb{P}^2$, and $\mathbb{P}^3$ approximations are computed up to $T=1$ on successively refined $N\times N$ uniform meshes using the four-entropy EPO-CDG method ($r\in\{0,0.50/\gamma,0.90/\gamma,0.97/\gamma\}$) and a CFL number of $0.01$. The resulting density errors and convergence rates are shown in \cref{tab:epocdg-2d-euler-density-wave-accuracy}. For $k=1,2,3$, the errors decrease rapidly under mesh refinement, attaining or exceeding the designed $(k+1)$th-order on finer meshes. These results confirm that enforcing multiple entropy constraints preserves the high-order accuracy of the scheme for smooth Euler flows.
\end{exmp}

\begin{table}[!htbp]
	\vspace{1em}
	\centering
	\setlength{\tabcolsep}{7pt}
	\renewcommand{\arraystretch}{1.04}
	\small
	\caption{Density errors and convergence rates for the 2D Euler density-wave problem in \cref{exmp:2d-euler-density-wave} at $T=1$, computed using the four-entropy $\mathbb{P}^k$-based EPO-CDG method.}
	\label{tab:epocdg-2d-euler-density-wave-accuracy}
	\begin{tabular}{@{}ccrcrcrc@{}}
		\toprule
		$k$ & $N$ &
		$L^1$ error & rate &
		$L^2$ error & rate &
		$L^\infty$ error & rate \\
		\midrule
		1 & 20 & $5.38\mathrm{e}{-2}$ & -- & $6.00\mathrm{e}{-2}$ & -- & $8.62\mathrm{e}{-2}$ & -- \\
		 & 40 & $8.45\mathrm{e}{-3}$ & 2.67 & $9.58\mathrm{e}{-3}$ & 2.65 & $1.47\mathrm{e}{-2}$ & 2.55 \\
		 & 80 & $1.03\mathrm{e}{-3}$ & 3.04 & $1.26\mathrm{e}{-3}$ & 2.93 & $2.19\mathrm{e}{-3}$ & 2.75 \\
		 & 160 & $1.27\mathrm{e}{-4}$ & 3.01 & $1.61\mathrm{e}{-4}$ & 2.96 & $2.94\mathrm{e}{-4}$ & 2.90 \\
		\midrule
		2 & 20 & $1.16\mathrm{e}{-3}$ & -- & $1.22\mathrm{e}{-3}$ & -- & $1.59\mathrm{e}{-3}$ & -- \\
		 & 40 & $7.56\mathrm{e}{-5}$ & 3.93 & $8.09\mathrm{e}{-5}$ & 3.92 & $1.31\mathrm{e}{-4}$ & 3.60 \\
		 & 80 & $4.70\mathrm{e}{-6}$ & 4.01 & $5.26\mathrm{e}{-6}$ & 3.94 & $9.34\mathrm{e}{-6}$ & 3.81 \\
		 & 160 & $3.01\mathrm{e}{-7}$ & 3.97 & $3.46\mathrm{e}{-7}$ & 3.93 & $7.14\mathrm{e}{-7}$ & 3.71 \\
		\midrule
		3 & 20 & $1.03\mathrm{e}{-4}$ & -- & $1.22\mathrm{e}{-4}$ & -- & $2.39\mathrm{e}{-4}$ & -- \\
		 & 40 & $4.63\mathrm{e}{-6}$ & 4.47 & $5.22\mathrm{e}{-6}$ & 4.54 & $9.23\mathrm{e}{-6}$ & 4.70 \\
		 & 80 & $3.57\mathrm{e}{-8}$ & 7.02 & $4.07\mathrm{e}{-8}$ & 7.00 & $8.71\mathrm{e}{-8}$ & 6.73 \\
		\bottomrule
	\end{tabular}
\end{table}

\begin{figure}[!htbp]
	\centering
	\begin{minipage}[t]{0.485\textwidth}
		\centering
		\includegraphics[width=\linewidth]
		{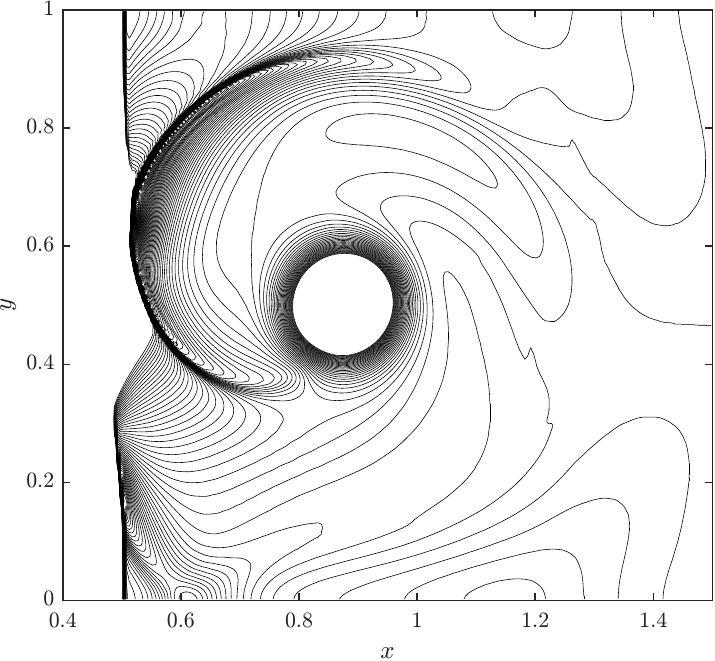}
	\end{minipage}\hfill
	\begin{minipage}[t]{0.485\textwidth}
		\centering
		\includegraphics[width=\linewidth]
		{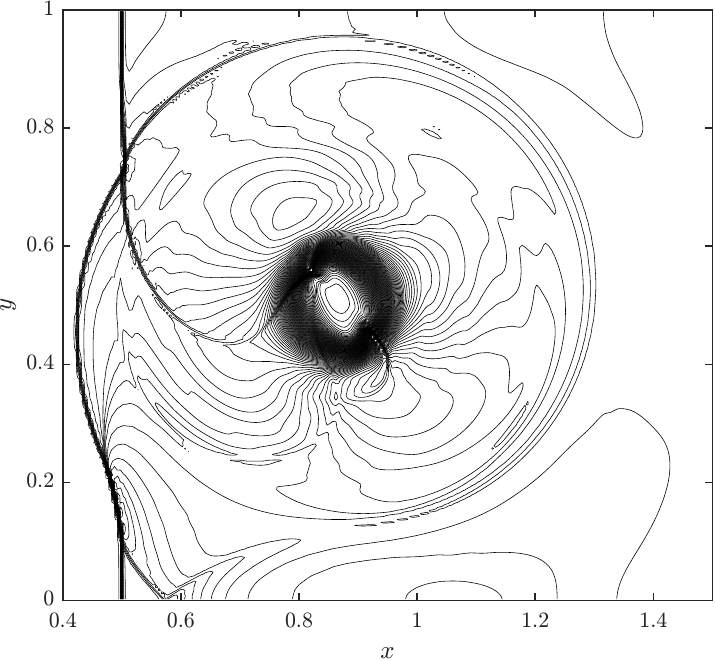}
	\end{minipage}
	\caption{Pressure contours at $t=0.529$ for the mild (left) and strong (right) shock-vortex interactions in \cref{exmp:2d-euler-shock-vortex}, using 91 contour lines over $[1.19,1.36]$ and $[0.02,1.39]$, respectively.}
	\label{fig:epocdg-2d-euler-shock-vortex}
\end{figure}

\begin{exmp}[Shock--vortex interaction]\label{exmp:2d-euler-shock-vortex}
		We consider the interaction of a Mach-$1.1$ stationary shock with an isentropic vortex on $\Omega = [0,2] \times [0,1]$. Initially, the stationary shock is located at $x = 0.5$, with the undisturbed upstream state $(\rho, u, v, p) = (1,\, 1.1\sqrt{\gamma},\, 0,\, 1)$ for $x < 0.5$ and the corresponding post-shock state for $x > 0.5$ determined by the Rankine--Hugoniot relations. An isentropic vortex centered at $(x_c, y_c) = (0.25, 0.5)$ is superimposed onto the upstream flow. Letting $\zeta := \sqrt{(x - x_c)^2 + (y - y_c)^2}/r_c$, $\Theta := p/\rho$, and $S := \log(p/\rho^\gamma)$, the vortex perturbations are given by
	$$
	(\delta u, \delta v)
	=
	\frac{\epsilon}{r_c}\mathrm{e}^{\alpha(1-\zeta^2)}
	(y - y_c,\, x_c - x),
	\qquad
	\delta \Theta
	=
	-\frac{(\gamma - 1)\epsilon^2}{4\alpha\gamma}
	\mathrm{e}^{2\alpha(1-\zeta^2)},
	\qquad
	\delta S = 0,
	$$
	which are added to the velocity, temperature, and specific entropy of the undisturbed upstream state, respectively. Inflow and outflow boundary conditions are imposed at the left and right boundaries, respectively, while reflecting boundary conditions are applied at the upper and lower boundaries. Both vortex tests are computed up to $T = 0.529$ using the $\mathbb{P}^3$ EPO-CDG method on a $400 \times 200$ uniform mesh, with four normalized Harten entropy pairs ($r \in \{0,\, 0.50/\gamma,\, 0.90/\gamma,\, 0.97/\gamma\}$) and third-order relaxation SSP Runge--Kutta startup followed by a third-order variable-step SSP multistep method. We simulate the following two vortices problem with different strength:
	\begin{itemize}[leftmargin=*]
		\item \textbf{Mild vortex}~\cite{PengWuYuan2025OECDG}.
		$(\epsilon,r_c,\alpha)=(0.3,0.05,0.204)$.
		
		\item \textbf{Strong vortex}~\cite{CuiDingWu2024OCAD}.
		$(\epsilon,r_c,\alpha)=(1.378106,0.04,0.204)$.
	\end{itemize}
	
	\Cref{fig:epocdg-2d-euler-shock-vortex} compares the flow structures for the mild and strong vortices. The mild vortex induces only a slight deformation of the shock front, whereas the strong vortex induces significant shock bending and rich secondary vortical structures. In both cases, EPO-CDG captures the complex wave interactions sharply without visible spurious oscillations.
\end{exmp}

For both astrophysical jet problems below, computations are performed on the upper half-domain ($y \ge 0$) in \texttt{long double} precision. Full-domain plots of $\log(\rho)$ are obtained by reflection across $y=0$. In PCDG, POCDG, and EPO-CDG, we evaluate the global discrete Harten entropy residuals $\mathscr R_r^{n+1}$ in \eqref{eq:global-discrete-entropy-residual} for $r\in\{0,0.50/\gamma,0.90/\gamma,0.97/\gamma\}$ after each variable-step SSP multistep update.

\begin{figure}[!htbp]
	\centering
	\includegraphics[width=0.99\textwidth]
	{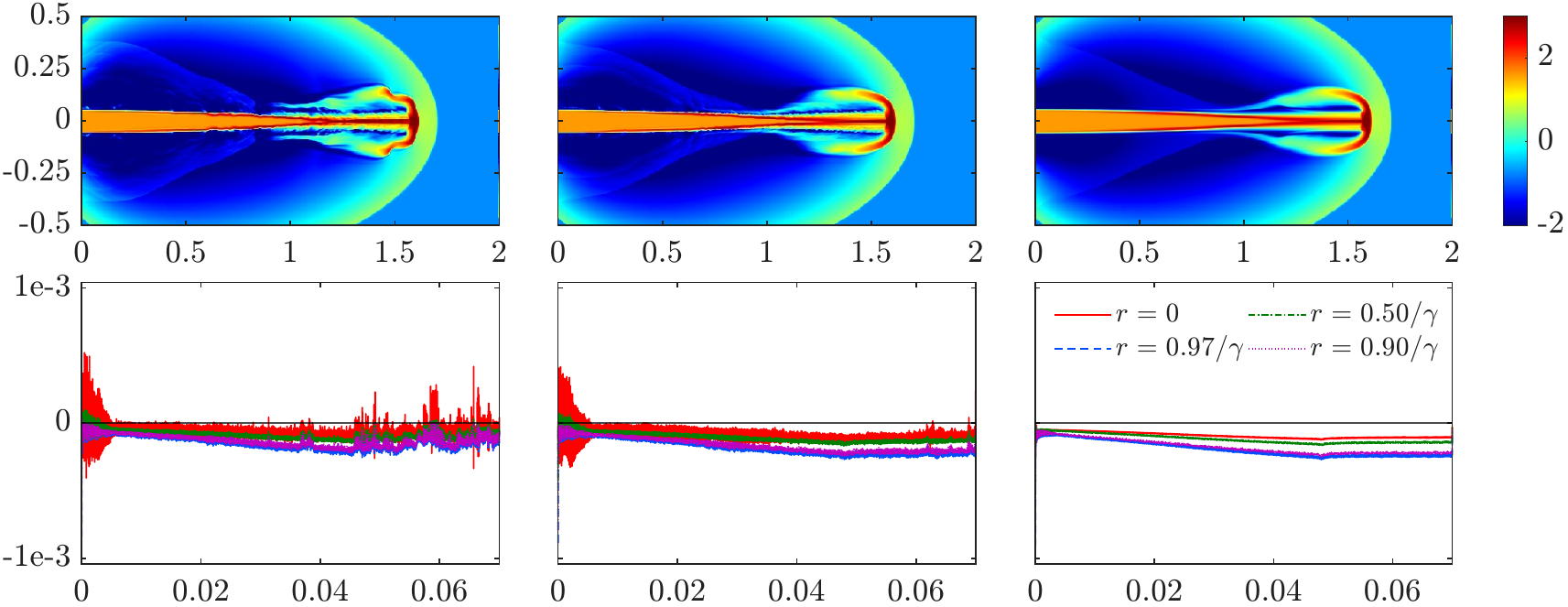}
	\caption{$\log(\rho)$ (top) and global discrete entropy residuals $\mathscr R_r^{n+1}$ (bottom) for the Mach 80 jet in \cref{exmp:2d-euler-mach80}. Columns from left to right show PCDG, POCDG, and EPO-CDG.}
	\label{fig:epocdg-2d-euler-mach80}
\end{figure}

\begin{exmp}[The Mach 80 jet problem \cite{ha2005numerical}]\label{exmp:2d-euler-mach80}
	On $\Omega=[0,2]\times[0,0.5]$, we set $\gamma=5/3$ and initialize the state as $(\rho,u,v,p)=(0.5,0,0,0.4127)$. At the left inflow boundary $x=0,\ 0\le y<0.05$, we prescribe the fixed state $(\rho,u,v,p)=(5,30,0,0.4127)$. A reflecting boundary condition is imposed along $y = 0$, while outflow conditions are prescribed at the upper and right boundaries as well as along $x = 0$ for $y \ge 0.05$. We use the $\mathbb{P}^2$ elements on a $640\times160$ mesh and compute up to $T=0.07$. As shown in \cref{fig:epocdg-2d-euler-mach80}, PCDG displays severe spurious oscillations near the jet core and shear layers. Although POCDG suppresses most of these spurious oscillations, unphysical structures persist in the same regions, accompanied by positive entropy residuals that violate the discrete entropy inequalities. In contrast, EPO-CDG resolves the jet bow shock and shear-layer structures cleanly while keeping all four entropy residuals nonpositive throughout the computation. 
\end{exmp}

\begin{figure}[!htbp]
	\centering
	\includegraphics[width=0.99\textwidth]
	{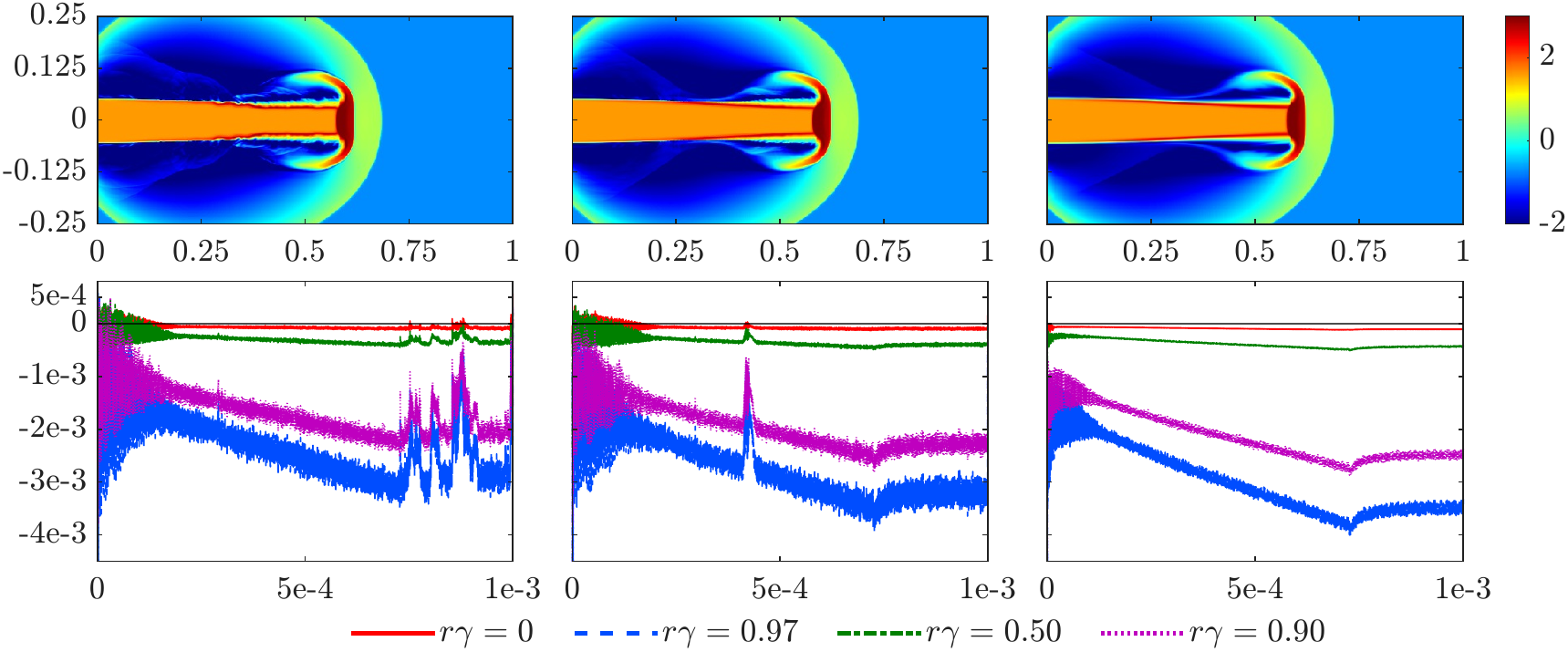}
	\caption{$\log(\rho)$ (top) and global discrete entropy residuals $\mathscr R_r^{n+1}$ (bottom) for the Mach 2000 jet in \cref{exmp:2d-euler-mach2000}. Columns from left to right show PCDG, POCDG, and EPO-CDG. }
	\label{fig:epocdg-2d-euler-mach2000}
\end{figure}

\begin{exmp}[Mach 2000 jet \cite{zhang2010positivity}]\label{exmp:2d-euler-mach2000}
	We use the same configuration as in the Mach 80 jet problem, except that the domain is set to $\Omega=[0,1]\times[0,0.25]$, the extreme inflow state is prescribed as $(\rho,u,v,p)_{\rm jet}=(5,800,0,0.4127)$, and the $\mathbb{P}^2$ approximations are computed up to $T = 0.001$ on a $640 \times 160$ uniform mesh. \Cref{fig:epocdg-2d-euler-mach2000} illustrates consistent behavior at this much higher Mach number. In \Cref{fig:epocdg-2d-euler-mach2000}, POCDG and EPO-CDG suppress the spurious oscillations observed near the bow shock in PCDG, and EPO-CDG resolves the shock and shear layers more clearly than the other two. Moreover, whereas some global entropy residuals for PCDG and POCDG become positive, all four residuals for EPO-CDG remain nonpositive at every time step.
\end{exmp}

\begin{figure}[!htbp]
    \centering
    \includegraphics[width=0.99\textwidth]{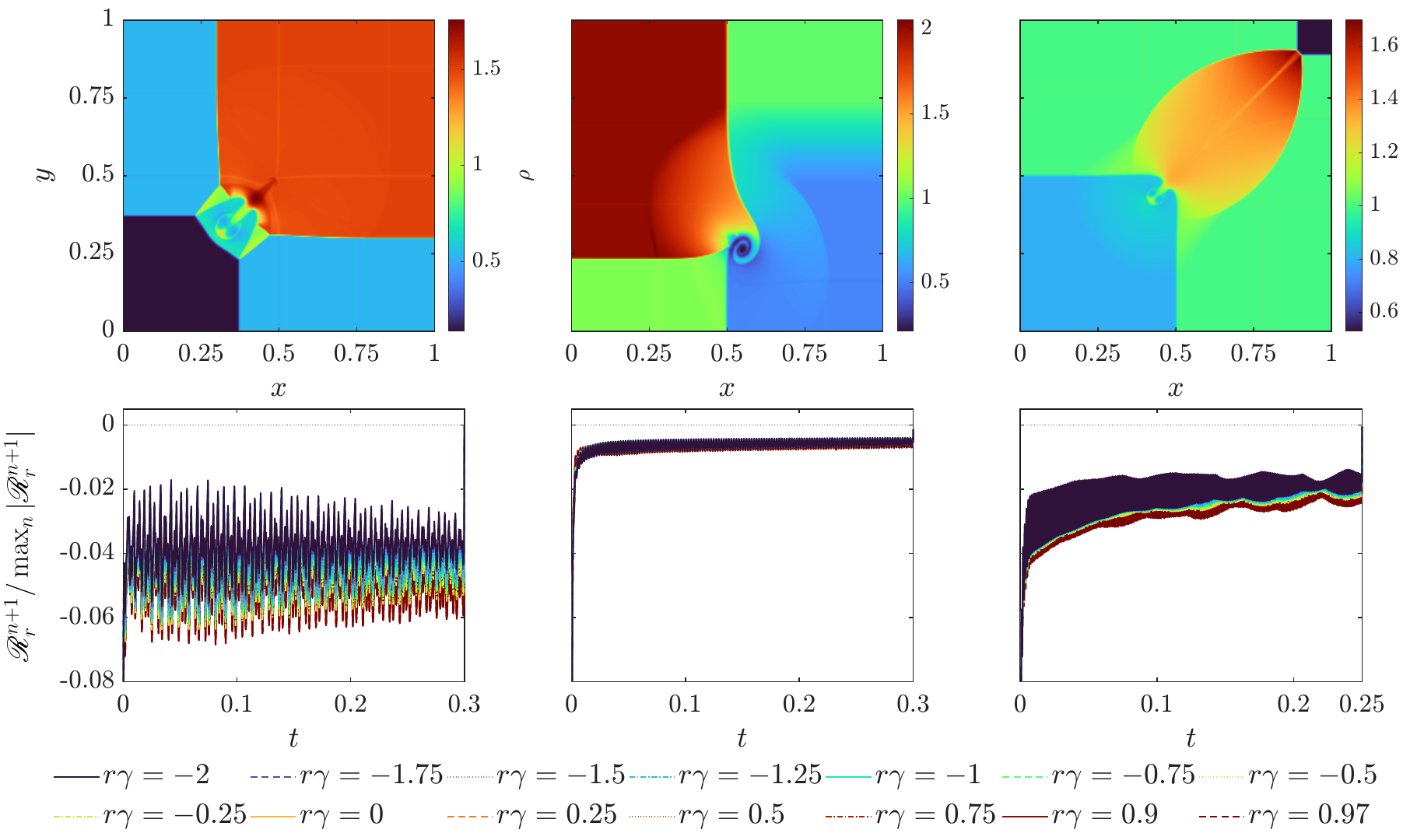}
    \caption{EPO-CDG density (top) and normalized fully discrete entropy residuals (bottom) for the three Riemann problems in \cref{exmp:2d-euler-riemann}. Columns from left to right correspond to cases (i)-(iii). }
    \label{fig:epocdg-2d-euler-riemann}
\end{figure}

\begin{exmp}[Riemann problems \cite{lax1998solution}]\label{exmp:2d-euler-riemann}
We consider three 2D Riemann problems on $\Omega=[0,1]^2$ with outflow boundary conditions on all boundaries. For each problem, the primitive-variable vector $\bm V=(\rho,u,v,p)$ is initially piecewise constant across the four quadrants defined by $x=0.5$ and $y=0.5$, denoted by $\mathrm{UL}$ (upper-left), $\mathrm{UR}$ (upper-right), $\mathrm{LL}$ (lower-left), and $\mathrm{LR}$ (lower-right). All simulations are performed using $\mathbb{P}^3$ EPO-CDG approximations on a $256 \times 256$ uniform mesh up to final times $T = 0.3$, $0.3$, and $0.25$ for configurations (i)--(iii), respectively, with 14 normalized Harten entropy constraints chosen such that $r\gamma \in \{-2 + 0.25j : j = 0, \ldots, 11\} \cup \{0.90, 0.97\}$. The initial states are given by:
	\begin{enumerate}[
		label=(\roman*),
		leftmargin=2em,
		labelsep=0.5em,
		nosep
		]
		\item {$\bm V_{\rm UL}=(0.5323,1.206,0,0.3)$, $\bm V_{\rm UR}=(1.5,0,0,1.5)$,\\
		$\bm V_{\rm LL}=(0.138,1.206,1.206,0.029)$, and $\bm V_{\rm LR}=(0.5323,0,1.206,0.3)$.}
		\item {$\bm V_{\rm UL}=(2,0,-0.3,1)$, $\bm V_{\rm UR}=(1,0,-0.4,1)$,\\
		$\bm V_{\rm LL}=(1.0625,0,0.2145,0.4)$, and $\bm V_{\rm LR}=(0.5197,0,-1.1259,0.4)$.}
		\item {$\bm V_{\rm UL}=(1,0.7276,0,1)$, $\bm V_{\rm UR}=(0.5313,0,0,0.4)$,\\
		$\bm V_{\rm LL}=(0.8,0,0,1)$, and $\bm V_{\rm LR}=(1,0,0.7276,1)$.}
	\end{enumerate}
	
	\Cref{fig:epocdg-2d-euler-riemann} shows that EPO-CDG cleanly resolves interacting shocks, contact discontinuities, and complex shear structures across all three test cases. For all three configurations, all 14 global discrete entropy residuals remain nonpositive throughout the time evolution.
\end{exmp}

\section{Conclusions} \label{sec:conclusions} 
	In this work, we developed the EPO-CDG framework, a high-order central discontinuous Galerkin method that combines fully discrete multi-entropy stability and strict bound preservation with oscillation suppression on staggered overlapping meshes. Under a unified CFL condition, the forward Euler (FE) update of the cell averages in CDG schemes satisfies both weak BP and weak discrete entropy inequalities for the prescribed family of convex entropy pairs.
	
	Conventional entropy-stable schemes typically rely on interface numerical fluxes and SBP flux differencing to achieve discrete entropy dissipation. Because CDG methods evolve solutions on staggered meshes without interface numerical fluxes, standard flux-differencing techniques do not directly apply, leaving high-order entropy stability a long-standing theoretical challenge for CDG discretizations. To overcome this difficulty, we established an intrinsic convex decomposition on each subcell that expresses the FE update of the cell averages as a convex combination of subcell quadrature values and directional LF intermediate states in CDG schemes. 
	By proving that these states are integral averages of self-similar Riemann solutions over the wave fan, we obtain the cell-average entropy inequalities from Jensen's inequality and the continuous entropy inequalities without interface numerical fluxes or Riemann solvers, while bound preservation follows from the convexity of the admissible set.
	To extend these properties from the cell averages to the high-order polynomial approximation, we constructed a weak-to-strong EPO lifting via local convex scaling. This lifting preserves the cell averages while enforcing the prescribed bounds at the designated quadrature points, satisfying multiple discrete entropy inequalities, and suppressing spurious Gibbs oscillations through scale-invariant primal-dual jump damping. Coupling this spatial lifting with variable-step-size strong-stability-preserving multistep time discretization yields high-order fully discrete schemes with provable multi-entropy stability and bound preservation.
	
	These results establish that provable fully discrete multi-entropy stability and bound preservation can be achieved on staggered overlapping meshes without summation-by-parts operators or interface numerical fluxes. Extensive numerical experiments confirm the optimal high-order accuracy of the proposed EPO-CDG framework for smooth solutions, as well as its robustness in extreme regimes. In particular, tests on nonconvex scalar conservation laws demonstrate that enforcing a single entropy inequality is often insufficient to suppress entropy-violating wave structures, whereas imposing multiple entropy constraints reliably recovers the physically correct composite waves. For the compressible Euler equations—including near-vacuum blast waves, strong shock interactions, and extreme Mach-2000 astrophysical jets—the method resolves complex wave structures sharply, strictly preserves the positivity of density and pressure, and maintains nonpositive global discrete entropy residuals for all prescribed entropies. In the further work, we will focus on optimizing subcell convex decompositions to relax the CFL restriction, developing systematic strategies for selecting entropy families in general hyperbolic systems, and extending the EPO-CDG framework to unstructured meshes, relativistic hydrodynamics, magnetohydrodynamics, and nonconservative hyperbolic systems.
	
	\appendix
	
\section{Admissibility and entropy dissipation of the LF intermediate state}
\label{app:weak-ES-BP}

This appendix proves that the directional LF intermediate state $\bm H_\ell$ belongs to $\mathcal G$ and satisfies the entropy inequalities stated in \Cref{lem:epocdg-H-admissibility}.
For the semidiscrete local Lax--Friedrichs scheme, Tadmor~\cite[Example~5.2]{tadmor2003entropy} established entropy stability under a numerical viscosity bounded by the maximum characteristic speed along an entropy-variable path. While this condition ensures entropy dissipation under implicit backward Euler time discretization, explicit time discretization requires sufficient spatial entropy dissipation to balance the temporal entropy production. In our setting, integrating a self-similar Riemann fan whose wave speeds are bounded by $\alpha_\ell$ and that satisfies the entropy condition for every $\eta^{(r)}\in\mathcal E$ directly yields the entropy inequality for the LF intermediate state. Under the explicit CFL condition, the updated cell average is a convex combination of admissible states. The convexity of the admissible set yields weak BP. The entropy estimate for the LF intermediate states and Jensen's inequality yield weak multi-entropy stability.

\begin{proof}[Proof of \Cref{lem:epocdg-H-admissibility}]
	Fix a coordinate direction $\ell$. Consider the 1D Riemann problem in that direction with left and right states $\bm U_L,\bm U_R$, flux $\bm f_\ell$, and wave-speed bound $\alpha_\ell$. Let $\bm R_\ell(\xi;\bm U_L,\bm U_R)$, with self-similar variable $\xi=x_\ell/t$, denote the corresponding entropy solution satisfying the entropy inequalities for every $\eta^{(r)}\in\mathcal E$. Since $\alpha_\ell$ bounds all waves in the Riemann fan, we have $\bm R_\ell(\xi)=\bm U_L$ for $\xi<-\alpha_\ell$ and $\bm R_\ell(\xi)=\bm U_R$ for $\xi>\alpha_\ell$. Integrating the conservation law $\partial_t \bm R_\ell + \partial_{x_\ell} \bm f_\ell(\bm R_\ell) = 0$ over the spacetime cone $\{(x_\ell,t):0<t<1,\ -\alpha_\ell t<x_\ell<\alpha_\ell t\}$ and applying the divergence theorem yields
	$$
	\frac{1}{2\alpha_\ell}
	\int_{-\alpha_\ell}^{\alpha_\ell}
	\bm R_\ell(\xi;\bm U_L,\bm U_R)\,d\xi
	=
	\frac12(\bm U_L+\bm U_R)
	-
	\frac{\bm f_\ell(\bm U_R)-\bm f_\ell(\bm U_L)}{2\alpha_\ell}
	=
	\bm H_\ell(\bm U_L,\bm U_R).
	$$
	Because $\mathcal G$ is a convex invariant domain and $\alpha_\ell$ bounds the complete Riemann fan, the Riemann-fan average on the left-hand side belongs to $\mathcal G$. Consequently, $\bm H_\ell(\bm U_L,\bm U_R)\in\mathcal G$. 
	For the entropy estimate, integrating the continuous entropy inequality $\partial_t\eta^{(r)}(\bm R_\ell) + \partial_{x_\ell}\mathcal Q_\ell^{(r)}(\bm R_\ell)\le0$ over the same spacetime cone and again applying the divergence theorem yields
	$$
	\frac{1}{2\alpha_\ell}
	\int_{-\alpha_\ell}^{\alpha_\ell}
	\eta^{(r)}(\bm R_\ell(\xi;\bm U_L,\bm U_R))\,d\xi
	\le
	\frac{\eta^{(r)}(\bm U_L)+\eta^{(r)}(\bm U_R)}{2}
	-
	\frac{\mathcal Q_\ell^{(r)}(\bm U_R)-\mathcal Q_\ell^{(r)}(\bm U_L)}{2\alpha_\ell}.
	$$
	Since $\eta^{(r)}$ is strictly convex, applying Jensen's inequality in conjunction with the averaging identity above gives
	$$
	\eta^{(r)}(\bm H_\ell(\bm U_L,\bm U_R))
	\le
	\frac{1}{2\alpha_\ell}
	\int_{-\alpha_\ell}^{\alpha_\ell}
	\eta^{(r)}(\bm R_\ell(\xi;\bm U_L,\bm U_R))\,d\xi .
	$$
	Combining these two inequalities completes the proof of \eqref{eq:epocdg-two-point-entropy-directional}.
\end{proof}

\section*{Supplement}
\textbf{The source code used to reproduce the numerical experiments in this paper is provided.}
\section*{Acknowledgements}
Google Gemini was used for English language editing. The authors reviewed the resulting edits and assume responsibility for all content. 
	\begin{changemargin}{-0.2cm}{-0.2cm}  %
		\renewcommand\baselinestretch{0.86}
		\bibliography{refs}
		\bibliographystyle{siamplain}
	\end{changemargin}
\end{document}

%% file: ex_shared.tex
\usepackage{lipsum}
\usepackage{amsfonts}
\usepackage{graphicx}
\usepackage{epstopdf}
\usepackage{algorithmic}
\ifpdf
\DeclareGraphicsExtensions{.eps,.pdf,.png,.jpg}
\else
\DeclareGraphicsExtensions{.eps}
\fi

\usepackage{booktabs}
\usepackage{hyperref}
\usepackage{stmaryrd}
\usepackage{bm}
\usepackage{caption}
\usepackage{subcaption}
\usepackage{multirow}
\usepackage{enumitem}

\usepackage{amsopn}

\allowdisplaybreaks

\usepackage{tikz, tikz-3dplot,forloop}

\usepackage[left=2.46cm,right=5.95cm,top=3.23cm,bottom=3.3cm,asymmetric]{geometry}

\newsiamthm{test}{Test}
\newsiamthm{assumption}{Assumption}

\newsiamremark{remark}{Remark}
\newsiamremark{hypothesis}{Hypothesis}
\crefname{hypothesis}{Hypothesis}{Hypotheses}
\newsiamthm{claim}{Claim}
\newsiamremark{exmp}{Example}

\headers{}{}

\title{A Unified High-Order Central Discontinuous Galerkin Framework with Fully Discrete Multi-Entropy Stability, Bound Preservation, and Oscillation Elimination
	\thanks{This work was partially supported by Science Challenge Project (No.~TZ2025007) and the Shenzhen Science and Technology Program (No.~JCYJ20250604144300001), the National Natural Science Foundation of China (No.~124B2022), and 
		the Postdoctoral Fellowship Program of CPSF (No.~2026T190820).}}

\author{Manting Peng \thanks{Department of Mathematics, Southern University of Science and Technology, Shenzhen, Guangdong 518055, China (\email{pengmt2024@mail.sustech.edu.cn}).} 
	\and Kailiang Wu \thanks{Corresponding author. Department of Mathematics and Shenzhen International Center for Mathematics, Southern University of Science and Technology, Shenzhen, Guangdong 518055, China (\email{wukl@sustech.edu.cn})}
		\and Jie Xu \thanks{ Department of Mathematics and Shenzhen International Center for Mathematics, Southern University of Science and Technology, Shenzhen, Guangdong 518055, China (\email{xuj9@sustech.edu.cn}).} 
 } 